\documentclass{article}

\usepackage[utf8]{inputenc}
\usepackage{amsthm,amsmath,amsfonts,amssymb,mathtools,xparse,bbm,graphicx,float,authblk}
\usepackage[margin=.9in]{geometry}
\usepackage{algpseudocode, algorithm}
\allowdisplaybreaks 

\usepackage{hyperref}
\hypersetup
{
	colorlinks,
	linkcolor = blue,
	citecolor = blue,
	urlcolor = blue,
	filecolor = blue,
	breaklinks,
}

\usepackage[
backend=bibtex,
style=numeric,
sorting=nyt,
maxbibnames=99,
]{biblatex}
\newtheorem{theorem}{Theorem}
\numberwithin{theorem}{section}
\newtheorem{proposition}[theorem]{Proposition}
\newtheorem{lemma}[theorem]{Lemma}
\newtheorem{corollary}[theorem]{Corollary}
\newtheorem{claim}[theorem]{Claim}
\theoremstyle{definition}
\newtheorem{definition}[theorem]{Definition}
\newtheorem{assumption}[theorem]{Assumption}
\theoremstyle{remark}
\newtheorem{remark}[theorem]{Remark}
\newtheorem{example}[theorem]{Example}

\newcommand{\RR}{\mathbb{R}}

\newcommand{\NN}{\mathbb{N}}

\newcommand{\set}[1]{{\left\{{#1}\right\}}}
\newcommand{\setc}[2]{\left\{\,#1 \;\vert\; #2\,\right\}}

\newcommand{\pow}[1]{\mathcal{P}(#1)}
\newcommand{\abs}[1]{\left| #1 \right|}
\newcommand{\norm}[1]{\left\lVert#1\right\rVert}
\newcommand{\linfnorm}[1]{\norm{#1}_{\infty}}

\newcommand{\PP}{\mathbb{P}}
\newcommand{\EE}{\mathbb{E}}
\newcommand{\EEbracket}[1]{\EE\left\{ #1 \right\}}
\newcommand{\Var}{\text{\rm Var}}
\newcommand{\Cov}{\text{\rm Cov}}
\newcommand{\ind}[1]{\mathbbm{1}\left(#1\right)}
\newcommand{\bin}[2]{\text{\rm Binomial}(#1, #2)}
\newcommand{\ber}[1]{\text{\rm Bernoulli}(#1)}

\newcommand{\normal}[2]{\mathcal{N}\left(#1, #2\right)}
\newcommand{\mvn}[2]{\mathrm{MVN}(#1, #2)}

\newcommand{\dks}[2]{d_\text{\rm KS}\left(#1, #2\right)}

\newcommand{\cplx}{\mathcal{L}}

\NewDocumentCommand{\gnp}{O{n} O{p}}{{\text {\bf G}}(#1, #2)}
\NewDocumentCommand{\xnp}{O{n} O{p}}{{\text {\bf X}}(#1, #2)}
\NewDocumentCommand{\multixnp}{O{n} O{p}}{{\text {\bf X}}(#1, \mathbf{#2})}
\NewDocumentCommand{\multixnpnotbold}{O{n} O{p}}{{\text {\bf X}}(#1, #2)}

\newcommand{\I}{\mathbb{I}}
\newcommand{\D}{\mathbb{D}}
\newcommand{\testzero}{\mathcal{H}_d}

\usepackage[usenames]{xcolor}
\newcommand{\tat}[1]{\textcolor{black}{{{}}#1}{}}

\newcommand{\gr}[1]{\textcolor{black}{{{}}#1}{}}

\title{\bf Goodness-of-fit via Count Statistics in \\Dense Random Simplicial Complexes}
\author[1]{Tadas Tem\v{c}inas}
\author[2]{Vidit Nanda}
\author[1]{Gesine Reinert}
\affil[1]{Department of Statistics, University of Oxford}
\affil[2]{Mathematical Institute, University of Oxford}

\begin{document}

\maketitle

\begin{abstract}
A key object of study in stochastic topology is a random simplicial complex. In this work we study a multi-parameter random simplicial complex model, where the probability of including a $k$-simplex, given the lower dimensional structure, is fixed. This leads to a conditionally independent probabilistic structure. This model includes the Erdős–Rényi random graph, the random clique complex as well as the Linial-Meshulam complex as special cases. The model is studied from both probabilistic and statistical points of view. We prove multivariate central limit theorems with bounds and known limiting covariance structure for the subcomplex counts and the number of critical simplices under a lexicographical acyclic partial matching. We use the CLTs to develop \gr{a} goodness-of-fit test for this random model and evaluate \gr{its} empirical performance. In order for the test to be applicable in practice, we also prove that the MLE estimators are asymptotically unbiased, consistent, uncorrelated and normally distributed.

\end{abstract}

\section{Introduction}
\gr{While complex data are often represented as graphs or networks \cite{newman2018networks}, t}here is a growing interest in modelling complex systems beyond pairwise interactions between nodes \cite{bick2021higher}. Simplicial complexes, \gr{often used in topological data analysis (TDA),} provide \gr{a rich mathematical representation of} higher-order networks \cite{bianconi}; other representations are detailed for example in \cite{barthelemy2022class}. As in network science, the statistical analysis of simplicial complex data \gr{relies on random complexes as null models.} In this paper, statistical understanding is achieved through a central limit theorem (CLT). Stochastic topology, which studies such random complexes, is partly motivated by its use in TDA \cite{bobrowski2018topology, bobrowski2022random}. Despite the growing popularity of TDA, random simplicial complexes are rarely used in practice, partly because of the lack of results in parameter estimation and asymptotic distribution of relevant statistics. This work aims to fill in the gap by providing \gr{a} probabilistic and statistical understanding of the so-called multiparameter random simplicial complex model $\multixnp$, which is rather general and includes several important special cases such as the Linial-Meshulam complex, the random clique complex $\xnp$ and the Bernoulli random graph $\gnp$. Here $\mathbf{p} = (p_1, p_2, \ldots , p_{n-1})$ is a vector of probabilities; $p_1$ is the probability of an edge, $p_2$ is the probability of a triangle to form a 2-simplex, \gr{and in general, $p_k$ is the probability of a $k$-vertex clique to form a $(k-1)-$simplex.} In particular  $\gnp$ is a special case of $\multixnp$ with  $\mathbf{p} = (p, 0, \ldots, 0).$ A formal definition is given in Section \ref{subsection:definitions}.

The motivation of this work is three-fold: 
\begin{enumerate}
    \item establishing results that allow the use of $\multixnp$ as a null-model in practical TDA applications;
    \item establishing an understanding of multivariate count statistics which have been instrumental in the study of algebraic invariants (such as Euler characteristic, homology groups, and the fundamental group) of random simplicial complexes; and finally,
    \item quantifying the extent to which discrete Morse theory can simplify homology computation in a random setting.
\end{enumerate}  

\gr{We address these} three points 
 in the \gr{so-called {\it dense} regime}, where the vector $\mathbf{p} \in (0,1)^{n-1}$ does not depend on $n$ \gr{(although some relaxation of this assumption is also discussed)}.
 \gr{A key observation which we establish in this paper is that 
  asymptotically the covariance matrix of all subcomplex counts in $\multixnp$ has rank one, which might be surprising. Even though the model has conditionally independent randomness in higher dimensions, this gets lost in the limit and so the behaviour of this random vector is analogous to a vector of subgraph counts in the $\gnp$ random graph model. Intuitively, this means that one subcomplex count asymptotically determines all the other ones. We shall also establish that the behaviour of the second count statistic - critical simplex counts of a lexicographical acyclic partial matching, as arising in discrete Morse theory - is slightly different.
 Next we cover the three motivations for this paper in more detail.}

\subsubsection*{Multiparameter complex as a null model}

In order to further statistical understanding of the model, we study maximum likelihood estimation in $\multixnp$ model as well as goodness-of-fit tests. We prove that the maximum likelihood estimator (MLE) for $\mathbf{p}$ is asymptotically unbiased and consistent. Also, we provide non-asymptotic bounds on the normal approximation error as well as the limiting covariance structure of the estimator. To our knowledge, this is the first study of the MLE for this random model. The classic abstract MLE results do not apply here because the observed simplex indicator variables are neither independent nor identically distributed\gr{; instead we employ results from \cite{temvcinas2021multivariate} which are based on a probabilistic technique called {\it Stein's method}}. 

\gr{Unfortunately, the vector of subcomplex counts  has limiting covariance matrix of rank one, which makes it an} unsuitable candidate\gr{s} for standard goodness-of-fit tests. Instead, we propose  critical simplex counts of a lexicographical acyclic partial matching \gr{as test statistic, because} this multivariate statistic has \gr{a nontrivial correlation} matrix \gr{also in the limit}. We analyse its performance empirically in a simulation study and show that it is able to distinguish between $\multixnp$ and a selection of geometric simplicial complex models. This is the first work \gr{proposing a} goodness-of-fit \gr{test} for \gr{the $\multixnp$ model.}

\subsubsection*{Distributional approximation of count statistics}

This work contains several probabilistic results, which provide some understanding of the probabilistic structure and which also motivate the choice of statistics for goodness-of-fit tests. In the study of random graphs, distributional approximations of subgraph counts play an important role \cite{rucinski1988small,barbour1989central}. Surprisingly, subcomplex counts in $\multixnp$, have not been studied from a distributional approximation point of view before, only from a large deviations point of view \cite{samorodnitsky2022large}. In this work we prove a multivariate CLT for subcomplex counts with non-asymptotic bounds on convergence rates with respect to both smooth test functions as well as convex set indicator functions. {Non-asymptotic bounds can be useful in applications since the observed data, however large, is never infinite.} \gr{Moreover we} prove a multivariate CLT with bounds on rates for \gr{critical simplex counts of a lexicographical acyclic partial matching} and derive a formula for the asymptotic covariance matrix, which generically has full rank, in contrast to the subcomplex counts. In order to prove a multivariate CLT for the subcomplex counts, we establish a multivariate CLT for generalised $U$-statistics under some assumptions, which might be of independent interest.

Proving the existence of certain subcomplexes have been a crucial step in showing results about homology, persistent homology, and the fundamental group in the regimes when these algebraic invariants are non-trivial with high probability \cite{ababneh2022maximal, roy2023random, newman2021one}. Distributional approximation of certain subcomplexes have also been used to deduce distributional approximation results for Betti number in the random clique complex \cite{kahle2013limit}. We hope that understanding the multivariate counts of arbitrary subcomplexes in $\multixnp$ as well as critical simplex counts under a lexicographical matching will pave a way to understand the multivariate distribution of Betti numbers in the dense regime.

\subsubsection*{Effectiveness of discrete Morse theory}

Discrete Morse theory \cite{forman2002user} provides a powerful, flexible and widely-used mechanism for simplifying the machine computation of simplicial homology and related algebraic invariants \cite{hmmn, mn, cgn, ripser, mtcog}. The basic idea is to construct a {\em partial matching} which pairs adjacent simplices subject to a global acyclicity constraint; the homology of the original complex is then entirely determined by (a chain complex constructed from) the unpaired simplices, \gr{the so-called  {\em critical} simplices}. Several strategies have been proposed for constructing acyclic partial matchings which admit relatively few number of critical simplices, thus greatly easing the linear-algebraic burden of computing homology. Simultaneously, considerable efforts have been invested in understanding the efficacy of such reductions in terms of {lowering} the number of simplices to consider --- the complexity of finding optimal matchings is discussed in \cite{joswig} while certain empirical calculations have been described in \cite{lutz}. In this paper we provide a multivariate normal approximation for vectors of critical simplex counts of {\em lexicographical} discrete Morse functions on $\xnp$. 
\gr{This acyclic partial matching is a standard tool in TDA which is used in some of the popular TDA software packages, for example  \cite{ripser}.} \gr{The multivariate normal approximation in our paper thus} furnishes a convenient benchmark for testing the efficacy of homology-preserving reductions via discrete Morse theory on random simplicial complexes.

\subsection*{Main results}

When presenting our main results here, we assume that the reader is familiar with random simplicial complexes and the multiparameter random simplicial complex $\multixnp$. {The subsequent Section \ref{section:top_preliminaries} provides the necessary background for the readers who are less familiar with these notions.}

There are two central theoretical contributions of this work: one regarding the properties of the MLE in $\multixnp$ and one regarding the distributional approximation of the mentioned count statistics. We first state a simplified version of our MLE results here. We state the results asymptotically for simplicity and clarity but we prove the results non-asymptotically with explicit constants. To see the results in full generality and detail, we refer to Theorems \ref{theorem:mle} and \ref{theorem:mle_normality} as well as Lemma \ref{lemma:mle_covariance}.

\begin{theorem}
    Consider the random simplicial complex $\multixnp$. Assume there is a known and fixed integer $d \geq 1$ such that $\mathbf{p} = (p_1, p_2, \ldots, p_{n-1})$ is a vector of probabilities where $p_i = 0$ for $i > d$ and $p_i \in (0,1)$ for $i \leq d$. Let $\hat{p}_i$ be the MLE of $p_i$ for some fixed $i \leq d$. Let $W_i \gr{= W_i{(n)}} = \binom{n}{i+1}^{\frac{1}{2}}(\hat{p}_i - p_i)$ and $Z_i \sim \normal{0}{(1-p_i)\prod_{j=1}^{i} p_j^{-\binom{i+1}{j+1}}}$. Then
    \begin{enumerate}
        \item \[ \lim_{n \to \infty} \EEbracket{\hat{p}_i} = p_i \text{ and }  \lim_{n \to \infty} \Var(\hat{p}_i) = 0; \]
        \item for any fixed $1 \leq i < j \leq d$ we have \[\lim_{n \to \infty}\Cov(W_i, W_j) = 0;\]
        \item for any fixed $1 \leq i \leq d$, \[ \lim_{n \to \infty} \sup_{x \in \RR} \abs{\PP(W_i \leq x) - \PP(Z_i \leq x)} = 0.\]
    \end{enumerate}
\end{theorem}
\gr{Here we recognize 
$\sup_{x \in \RR} \abs{\PP(W_i \leq x) - \PP(Z_i \leq x)}
=\sup_{x \in \RR} \abs{\EE \{ \mathbbm{I} (W_i \in (-\infty, x] )\}  -
 - \EE \{ \mathbbm{I} (Z_i \in (-\infty, x] )\} } $ as the Kolmogorov-distance between the distribution of $W_i$ and the distribution of $Z_i$, which is based on indicator test functions $\mathbbm{I} (z \in (-\infty, x] )$, taking the value 1 if $z \in (-\infty, x]$ and 0 otherwise. A natural multivariate generalisation of this distance is to take the supremum over indicator functions of convex sets. Similarly, taking instead the supremum over Lipschitz(1)-test functions yields Wasserstein distance.}

\gr{Next} we state a theorem that contains a simplified version of the multivariate normal approximation results, \gr{for both, simplex counts and and critical simplex counts}. \tat{Even though the two random vectors are related from a topological point of view, a CLT for one of them does not automatically translate into a CLT for the other: the dependence structure in the two cases are quite different, requiring slightly different proofs.} We quantify the approximation error in terms of both smooth test functions and convex set indicator functions. Even though an error vanishing asymptotically implies convergence in distribution in both cases, for finite $n$, the convex set indicator functions give us a stronger result that is also of interest in practice when, for example, estimating confidence regions. For the non-simplified version of the results we refer the reader to Theorems \ref{theorem:subcomplex_counts_limiting_cov}, \ref{theorem:crit_lexi_limiting_cov}, \ref{theorem:crit_lexi_approx},  Corollary \ref{theorem:subcomplex_counts}, \gr{and Lemma \ref{theorem:crit_lexi_limiting_cov}}. \gr{For the statement we l}et $k'$ be the largest number such that for all $i \leq k'$ we have $p_i = 1$; $k' = 0$  if none of the probabilities \gr{equal 1}.

\begin{theorem}\label{theorem:main_CLT}
    Consider the random simplicial complex $\multixnp$.  Let $W^{(1)} \in \RR^d$ be an appropriately scaled and centered vector of subcomplex counts. Let $W^{(2)} \in \RR^d$ be an appropriately scaled and centered count vector of simplices that are critical under a lexicographical acyclic partial matching. Let $Z \sim \mvn{0}{\text{\rm Id}_{d \times d}}$ and $\Sigma_i$ be the covariance matrix of $W^{(i)}$ for each $i$. Let $h: \RR^d \to \RR$ be three times partially differentiable function whose third partial derivatives are Lipschitz continuous and bounded. Let $\mathcal{K}$ be the class of convex sets in $\RR^d$.
    \begin{enumerate}
        \item There is an explicit constant $B_{\ref{theorem:subcomplex_counts}} > 0$, independent of $n$ such that \[
        \abs{\EE h(W^{(1)}) - \EE h(\Sigma_1^{\frac{1}{2}}Z)} \leq \sup_{i, j, k \in [d]} \linfnorm{\frac{\partial^3 h}{\partial x_i \partial x_j \partial x_k}}  B_{\ref{theorem:subcomplex_counts}}  n^{-\frac{k'+2}{2}}
        \] and 
         \[ \sup _{A \in \mathcal{K}}|\PP(W^{(1)} \in A)-\PP(\Sigma_1^{\frac{1}{2}}Z \in A)| \leq 2^{\frac{7}{2}} 3^{-\frac{3}{4}}d^{\frac{3}{16}}B_{\ref{theorem:subcomplex_counts}}^{\frac{1}{4}} n^{-\frac{k'+2}{8}}.\]
         Moreover, for any $1 \leq i < j \leq d$ we have \[\lim_{n \to \infty} (\Sigma_1)_{i,j} = 1.\]

        \item There is an implicit constant $B_{\ref{theorem:crit_lexi_approx}.1} >0$ independent of $n$ and a natural number $N_{\ref{theorem:crit_lexi_approx}.1}$ such that for any $n \geq N_{\ref{theorem:crit_lexi_approx}.1}$ we have \[\abs{\EE h(W^{(2)}) - \EE h(\Sigma_2^{\frac{1}{2}}Z)} \leq B_{\ref{theorem:crit_lexi_approx}.1} \sup_{i, j, k \in [d]} \linfnorm{\frac{\partial^3 h}{\partial x_i \partial x_j \partial x_k}} n^{-\frac{k'+2}{2}}.\]
         
        Also, there is an implicit constant $B_{\ref{theorem:crit_lexi_approx}.2} >0$ independent of $n$ and a natural number $N_{\ref{theorem:crit_lexi_approx}.2}$ such that for any $n \geq N_{\ref{theorem:crit_lexi_approx}.2}$ we have\[\sup _{A \in \mathcal{K}}|\PP(W^{(2)} \in A) -\PP(\Sigma_2^{\frac{1}{2}}Z \in A)| \leq B_{\ref{theorem:crit_lexi_approx}.2}n^{-\frac{k'+2}{8}}.\]

        Moreover, for any $1 \leq i < j \leq d$ there is an explicit function $\sigma_{\infty}(i,j)$ such that \[\lim_{n \to \infty} (\Sigma_2)_{i,j} = \sigma_{\infty}(i,j).\]
        \gr{If $p_k \in (0,1)$ for all $k > k'$, then then $0 \leq \sigma_{\infty}(i,j) < 1.$}
    \end{enumerate}
\end{theorem}
We also propose a testing procedure based on critical simplex counts and empirically evaluate its performance. The procedure is described in Algorithm \ref{algorithm:test} and the results are presented in Section \ref{subsection:testing}.

\subsection*{Related work}

Although subcomplex counts in $\multixnp$ are natural objects to consider, surprisingly, they have not yet been extensively studied. Nonetheless, a recent paper studies large deviation bounds for subcomplex counts in the $\multixnp$ model \cite{samorodnitsky2022large}. {Another recent paper \cite{eichelsbacher2023simplified} proves a CLT for a special kind of subcomplex counts in a special case of the $\multixnp$ model, the number of isolated faces in the Linial-Meshulam random simplicial complex.} \gr{As s}ubcomplex counts can be seen as an example of generalised $U$-statistics,  \gr{r}esults from \cite[Section 10.7]{janson1991asymptotic} and \cite[Chapter 11.3]{janson1997gaussian} could, with some work, in principle be adapted to prove convergence in distribution for the multivariate subcomplex counts; \gr{in} \cite{de1996central} the author extends the results from \cite{janson1991asymptotic} to study generalised $U$-statistics not only for graphs but also for hypergraphs {and therefore for simplicial complexes in particular}; , \gr{but no non-asymptotic bounds on the approximation error are available}. Recently Zhang \cite{zhang2022berry} used the exchangeable pair approach in Stein's method to derive a univariate CLT for the graph version of generalised $U$-statistics. These results only {apply to} graph statistics rather than the higher order type of statistics that we are interested in here.

Counts of simplices that are critical under a lexicographical acyclic partial matching have been studied for the random clique complex $\xnp$ \cite[Section 4]{temvcinas2021multivariate}. For the $\multixnp$ model, only the expected value has been studied so far \cite[Section 8]{bauer2021parameterized}.

There has been some limited work on statistics for random simplicial complexes. In \cite{bobrowski2022random, zuev2015exponential} the authors have noted a sufficient statistic for the $\multixnp$ model but the MLE had not been investigated. Goodness-of-fit tests using count statistics have been proposed in the context of random graphs \cite{kaur2020higher, bubeck2016testing, liu2021phase}. Topological goodness-of-fit tests have been suggested for  point processes \cite{biscio2020testing}, sliced spatial data \cite{cipriani2023topology}, and cylindrical networks \cite{krebs2022functional}, \gr{but not, to the best of our knowledge, for} $\multixnp$.

\subsection*{Organisation and notation}

In Section \ref{section:top_preliminaries} we introduce topological definitions and the random models of interest. In Section \ref{section:stats} we study the properties of the MLE of the parameter $\mathbf{p}$ in  $\multixnp$ and perform a simulation study to verify the validity of the proposed goodness-of-fit test. In Section \ref{section:probabilistic_tools} we introduce and prove probabilistic results that are needed to prove Theorem \ref{theorem:main_CLT}. It also contains a multivariate CLT for generalised $U$-statistics that might be of independent interest. In Section \ref{section:subcomplex_counts} we study multivariate subcomplex counts and prove a multivariate CLT \gr{with an explicit bound; we also} describe the limiting covariance. In Section \ref{section:crit_lexi} we study critical simplex counts under a lexicographical acyclic partial matching. A multivariate CLT with \gr{an explicit bound} is proved and the limiting covariance described. Finally, the Appendix \ref{section:appendix} contains some technical computations that are needed for the proofs of Section \ref{section:crit_lexi}.

Throughout this paper we use the following notation. Given positive integers $n,m$ we write $[m,n]$ for the set $\set{m,m+1,\ldots,n}$ and $[n]$ for the set $[1,n]$. Given a set $X$ we write $\abs{X}$ for its cardinality, $\pow{X}$ for its powerset, and given a positive integer $k$ we write $\binom{X}{k} = \setc{t \in \pow{X}}{|t| = k}$ for the collection of subsets of $X$ which are of size $k$. For a function $f: \RR^d \to \RR$ we write  $\partial_{ij}f = \frac{\partial^2f}{\partial x_i \partial x_j}$ and $\partial_{ijk}f = \frac{\partial^3f}{\partial x_i \partial x_j \partial x_k}$. Also, we write $\abs{f}_k = \sup_{i_1, i_2, \ldots, i_k \in [d]} \linfnorm{\partial_{i_1 i_2 \ldots i_k}f}$ for any integer $k \geq 1$, as long as the quantities exist. Here $|| \cdot ||_\infty$ denotes the supremum norm. For a positive integer $d$ we define a class of test functions $h: \RR^d \to \RR$, as follows. We say $h \in \testzero$ iff $h$ is three times partially differentiable with third partial derivatives being Lipschitz and $\abs{h}_3 < \infty$. \gr{The} class of convex sets in $\mathbb{R}^d$ \gr{is denoted by $\mathcal{K}_d$.  W}hen the dimension of the space is clear from the context we \gr{may suppress the subscript $d$ for these sets}. The notation $\text{\rm Id}_{d \times d}$ denotes the $d \times d$ identity matrix. The vertex set of all graphs and simplicial complexes is assumed to be $[n]$, unless stated otherwise. We also use Bachmann-Landau asymptotic notation: we say $f(n) = O(g(n))$ iff $\limsup_{n \to \infty} \frac{|f(n)|}{g(n)} < \infty$, $f(n) = \Omega(g(n))$ iff $\liminf_{n \to \infty} \frac{f(n)}{g(n)} > 0$, and $f(n) = \omega(g(n))$ iff $\lim_{n \to \infty} \frac{f(n)}{g(n)} = \infty$. Given two real random variables $X, Y$ we write $\dks{X}{Y} \coloneqq \sup_{x \in \RR} \abs{\PP(X \leq x) - \PP(Y \leq x)}$ for the Kolmogorov-Smirnov distance \gr{between their distributions}.
Moreover, we use the following bounds on binomial coefficients
\begin{equation}\left(\frac{n}{k} \right)^k \le {n \choose k} \le \frac{n^k}{k!} < \left(\frac{ne}{k} \right)^k.
  \label{eq:binom}  
\end{equation}

\section{Topological preliminaries}\label{section:top_preliminaries}

For completeness, we include some topological preliminaries, see also \cite{temvcinas2021multivariate}. 

	\subsection{First definitions}\label{subsection:definitions}
	
	Firstly, we recall the notion of a simplicial complex \cite[Ch 3.1]{spanier}; these are special kinds of hypergraphs that provide  higher-dimensional generalisations of a graph and constitute data structures of interest across algebraic, applied, and computational topology.
	
	A simplicial complex $\cplx$ on a vertex set $[n]$ is a set of nonempty subsets of $[n]$ (i.e. $\varnothing \notin \cplx \subseteq \pow{[n]}$) such that the following properties are satisfied:
	\begin{enumerate}
		\item for each $v \in [n]$ the singleton $\set{v}$ lies in $\cplx$, and
		\item if $t \in \cplx$ and $s \subseteq t$ then $s \in \cplx$.
	\end{enumerate}
	
	The \textbf{dimension} of a simplicial complex $\cplx$ is $\dim(\cplx) \coloneqq \max_{s \in \cplx}|s| - 1$. Elements of a simplicial complex are called \textbf{simplices}. If $s$ is a simplex, then its dimension is $\dim(s) \coloneqq |s| - 1$. A simplex of dimension $k$ can be called a $k$-simplex. The $k$-skeleton of a simplicial complex $\cplx$ is the subcomplex $\cplx^{(k)} = \setc{s \in \cplx}{|s| \leq k + 1}$ of dimension $k$ as long as $k < \dim(\cplx)$. Note that the notion of one-dimensional simplicial complex is equivalent to the notion of a graph, with the vertex set $[n]$ and edges as one-dimensional simplices.
	
	{Given $k_1, k_2 \in \NN$ and two simplicial complexes $\cplx_1$ on $[k_1]$ and $\cplx_2$ on $[k_2]$, a \textbf{simplicial map} from $\cplx_1$ to $\cplx_2$ is a function $f: [k_1] \to [k_2]$ such that for any $s \in \cplx_1$ we have that $f(s) \in \cplx_2$. If $f$ is bijective, it is called a \textbf{simplicial isomorphism}. For any finite simplicial complex $\cplx$ on $[n]$ we denote the set of simplicial complexes on $[n]$ that are isomorphic to $\cplx$ by $[\cplx]$. Note that $[\cplx]$ can be smaller than the number of automorphisms on $\cplx$. For example, if $\cplx = \set{\set{1}, \set{2}, \set{1,2}}$ is just an edge then $[\cplx] = \set{\cplx}$. Even though the edge has two automorphisms (the identity and the one swapping the two vertices), both automorphisms give rise to the same simplicial complex. Indeed the swapping automorphism gives $\set{\set{2}, \set{1}, \set{2,1}}$, which is the same set as $\cplx$ itself.}
 
 {Assuming $k_1 \leq k_2$ it is also meaningful to count the number of copies of $\cplx_1$ in $\cplx_2$, which we define as the number of order-preserving injective simplicial maps from any element of $[\cplx_1]$ to $\cplx_2$. We call this the \textbf{subcomplex count} of $\cplx_1$ in $\cplx_2$. Subcomplex counts are investigated in Section \ref{section:subcomplex_counts} and concrete examples are given there.}
	
	\begin{definition}\label{definition:simplex_hollow_simplex}
	The $k$-simplex is equal to the set $\pow{[k+1]} \setminus \set{\varnothing}$, seen as a simplicial complex. That is,  the vertices are $\set{1,2,\ldots,k+1}$ and any subset of the vertex set is a simplex. The \textbf{hollow} $k$-simplex is equal to the set $\pow{[k+1]} \setminus \set{[k+1], \varnothing}$ seen as a simplicial complex. That is, a hollow $k$-simplex is the $(k-1)$-skeleton of the $k$-simplex and so its dimension is $k-1$.
	\end{definition}
	
	\begin{definition}\label{definition:costa-farber}
	Let $[n]$ be the vertex set. Let $\mathbf{p} = (p_1, p_2, \ldots, p_{n-1}) \in [0,1]^{n-1}$ be a finite sequence of parameters. We first build a random hypergraph $H$ where any subset $s \subseteq [n]$ of size $k$ is included independently of any other subsets with probability $p_{k-1}$. The \textbf{multi-parameter random simplicial complex} $\multixnp$ is the largest simplicial complex contained in $H$. That is, for any subset $s \subseteq [n]$, we have $s \in \multixnp$ iff $s \in H$ and for any $t \subseteq s$ we have $t \in H$. 
	\end{definition}
	
	\gr{The model} $\multixnp$ is a general random simplicial complex model, defined by Kahle in \cite{kahle2014topology} and extensively studied from a {topological rather than combinatorial point of view by}, amongst others, Costa and Farber in \cite{costa2016large, costa2017large}; it includes many families of random simplicial complexes as special cases. For example, if $\mathbf{p} = (p, 0, \ldots, 0)$, then we recover the $\gnp$ random graph. If we set  $\mathbf{p} = (p, 1, \ldots, 1)$, then we recover the random clique complex of $\gnp$, that has been extensively studied, for example, in \cite{kahle2009topology,kahle2014sharp,kahle2013limit}. If we set $\mathbf{p} = (1, \ldots 1, p_{k}, 0, \ldots, 0)$, then we recover the $k$-dimensional Linial-Meshulam random complex that has also been extensively studied before \cite{linial2006homological, linial2016phase}. For a survey of stochastic topology, see \cite{bobrowski2018topology, bobrowski2022random}.
	
	\subsection{Discrete Morse theory}

	A \textbf{partial matching} on a simplicial complex $\cplx$ is a collection
	\[
	V = \setc{(s, t)}{s \subseteq t \in \cplx \text{ and } |t| - |s| = 1}
	\] such that every simplex appears in at most one pair of $V$. A \textbf{$\mathbf{V}$-path} (of length $k \geq 1)$ is a sequence of distinct simplices of $\cplx$ of the following form:
	\[(s_1 \subseteq t_1 \supseteq s_2 \subseteq t_2 \supseteq \ldots \supseteq s_k \subseteq t_k)\]
	such that $(s_i, t_i) \in V$ and $|t_i| - |s_{i+1}| = 1$ for all $i \in [k]$. A $V$-path is called a \textbf{gradient path} if $k=1$ or $s_1$ is not a subset of $t_k$. A partial matching $V$ on $\cplx$ is called \textbf{acyclic} iff every $V$-path is a gradient path. Given a partial matching $V$ on $\cplx$, we say that a simplex $t \in \cplx$ is \textbf{critical} iff $t$ does not appear in any pair of $V$.

	For a one-dimensional simplicial complex, viewed as a graph, a partial matching $V$ is comprised of elements $( v; \{u,v\})$ with $v$ a vertex and $\{u, v\}$ an edge. A $V-$path is then a sequence of distinct vertices and edges 
	\[
	v_1, \set{v_1,v_2}, v_2, \set{v_2,v_3}, \ldots, v_k, \set{v_k,v_{k+1}}
	\]
	where each consecutive pair of the form $(v_i,\set{v_i,v_{i+1}})$ is constrained to lie in $V$.
	
	\medskip
	We refer the interested reader to \cite{forman2002user} for an introduction to discrete Morse theory and to \cite{mn} for \gr{an illustration of its use for simplifying} computation\gr{s} of persistent homology. \gr{This work addresses} how much computational improvement one should expect to gain on a random input when using a specific type of acyclic partial matching, defined below. 
	
	\begin{definition}\label{def:lexi_matching}
		Let $\cplx$ be  a simplicial complex and assume that the vertices are ordered by $[n] = \{1,\ldots,n\}$.
		For each simplex $s \in \cplx$ define
		\[
		I_{\cplx}(s) \coloneqq \{j \in [n] \mid j < \min(s) \text{ and } s \cup \{j\} \in \cplx\}. 
		\]
		Now consider the pairings
		\[
		s \leftrightarrow s \cup \{i\},
		\]
		where $i = \min I_{\cplx}(s)$ is the smallest element in the set $ I_{\cplx}(s)$, defined whenever $I_{\cplx}(s) \neq \varnothing$. We call this the {\bf lexicographical matching}.
	\end{definition}
	
	Due to the  $\min I_{\cplx}(s)$ construction in the lexicographical matching, the indices are decreasing along any path and hence \gr{the simplices form} a gradient path, showing that the lexicographical matching is indeed an acyclic partial matching on $\cplx$.

	\begin{figure}
		\centering
		\includegraphics[width=0.25\textwidth]{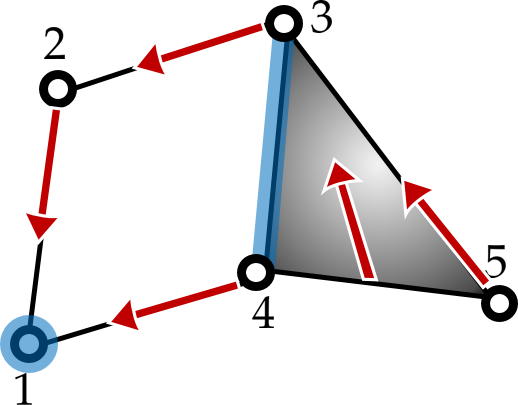}
		\caption{Lexicographical matching given by the red arrows. Critical simplices are highlighted in blue.}
		\label{fig:lexmorse}
	\end{figure}
	
	\begin{example}
		Consider the simplicial complex $\cplx$ depicted in Figure \ref{fig:lexmorse}. 
		The complex has 5 vertices, 6 edges and one two-dimensional simplex that is shaded in grey. The red arrows show the lexicographical matching on this simplicial complex: there is an arrow from a simplex $s$ to $t$ iff the pair $(s, t)$ is part of the matching. More explicitly, the lexicographical matching on $\cplx$ is \[V = \set{(\set{2}, \set{1,2}), (\set{3}, \set{2,3}), (\set{4}, \set{1,4}), (\set{5}, \set{3,5}), (\set{4, 5}, \set{3,4, 5})}.\] 
		Note that $\set{3,4}$ cannot be matched because the set $I_{\cplx}(\set{3,4})$ is empty. Also, in any lexicographical matching $\set{1}$ is always critical as there are no vertices with a smaller label and hence the set $I_{\cplx}(\set{1})$ is empty. So under this matching there are two critical simplices: $\set{1}$ and $\set{3,4}$, highlighted in blue in the figure. Hence, if we were computing the homology of this complex, considering only two simplices would be sufficient instead of all 12 which are in $\cplx$ -- a significant reduction in computational effort.
	\end{example}

\section{Maximum likelihood estimation and goodness-of-fit}\label{section:stats}

\subsection{Maximum likelihood estimation}

Recall the random simplicial complex model $\multixnp$ as given in Definition \ref{definition:costa-farber}. Here we study a maximum likelihood estimator under this random model assuming that there is some known $d$ such that $p_i \in (0,1)$ if $i \leq d$ and $p_i = 0$ if $i > d$. Therefore, we have $d$ parameters in this model to estimate, namely, $(p_1, p_2, \ldots, p_d) \in (0,1)^d$. With slight abuse of notation, we write $\mathbf{p} = (p_1, p_2, \ldots, p_d)$. Let $\cplx$ be a fixed simplicial complex on the vertex set $[n]$, which we consider to be our observed data. Recall Definition \ref{definition:simplex_hollow_simplex}, detailing the notions of a simplex and a hollow simplex. For each $i \in [n]$, let $s_i$ be the number of $i$-simplices and $h_i$ be the number of hollow $i$-simplices. Note that $s_i \leq h_i$ because every $i$-simplex contains a hollow $i$-simplex. Also note that $h_1 = \binom{n}{2}$ is just the number of pairs of vertices and hence it is a deterministic value.

\gr{
As it has been noted in \cite{bobrowski2022random, zuev2015exponential}, for any $i \in [d]$ a sufficient statistic for $p_i$ is the pair $(s_i, h_i)$. We therefore can estimate $p_i$ only for $i \leq I(\cplx)$. Note that if $h_i = 0$, then $h_j = 0$ and $s_j = 0$ for all $j \geq i$ and if $h_i \neq 0$, then $h_j \neq 0$ for all $j \leq i$.}
Write $I(\cplx) = \max \setc{i \in [d]}{h_i \neq 0}$;  in the dense regime, if  $d$ is not too large, we can prove that $I(\cplx) = d$ with high probability (see Lemma \ref{lemma:hollow_simplices_concentration}). Assuming $\mathbf{p} \in (0,1)^{d}$, the log-likelihood is well-defined and equal to 
\[l(\mathbf{p};\cplx) = \sum_{i=1}^{I(\cplx)} \left\{ s_i \ln(p_i) + (h_i - s_i) \ln(1 - p_i) \right\}.\]

For any $i \in [I(\cplx)]$ we have $\frac{\partial l}{\partial p_i} = \frac{s_i}{p_i} - \frac{h_i-s_i}{1-p_i}$. Solving the partials equal to zero and considering the boundary cases when $s_i=0$ or $s_i=h_i$, we get $\hat{p}_i = \frac{s_i}{h_i}$ for $i \leq I(\cplx)$. If $i > I(\cplx)$, then $s_i = 0$ and so the log-likelihood does not depend on $p_i$. In this case $\hat{p}_i$ is not unique; we choose to set $\hat{p}_i = 0$.  
\gr{Summarising,} we have
{\begin{align}\label{equation:mle}
    \hat{p}_i = \ind{h_i \neq 0} \frac{s_i}{h_i}.
\end{align}} 
For any $i \neq j \in [I(\cplx)]$ we have:
\begin{align*}
    \frac{\partial^2 l}{\partial p_i \partial p_j}(\hat{\mathbf{p}}) = 0; && \frac{\partial^2 l}{\partial p_i^2}(\hat{\mathbf{p}}) = -h_i^2\left( \frac{1}{s_i} + \frac{1}{h_i-s_i} \right) < 0.
\end{align*}

For $i \neq j \in [I(\cplx) + 1, d]$ the corresponding second partial \gr{derivatives all vanish; } the $I(\cplx) \times I(\cplx)$ Hessian is negative-definite, and hence $\hat{\mathbf{p}}$ indeed is a local (and global) maximum of the likelihood. Classical theorems from the theory of MLE do not apply here because \gr{the} simplex indicators are neither independent nor identically distributed. {Hence we derive asymptotic results \gr{for the MLE} in this paper.

For each $i \in [d]$ let $H_i$ be the random variable counting the number of hollow $i$-simplices and let $S_i$ be the random variable counting the number of solid $i$-simplices. Let $h_i$ be an integer such that $1 \leq h_i \leq \binom{n}{i+1}$. Then by definition of $\multixnp$  for any $i \in [d]$ we have that $S_i | \set{H_i = h_i} \sim \bin{h_i}{p_i}$. Using this property of $S_i$ being conditionally binomial we can prove that the estimator is asymptotically unbiased and consistent.  We start the investigation by proving concentration of $H_i$ for $i \in [d]$ and hence that $I(\cplx) = d$ with high probability.

For the rest of the section we write $p^{*} \coloneqq \min \setc{p_j}{j \in [d]}$ and $P_i = \prod_{j=1}^{i-1} p_j^{\binom{i+1}{j+1}}$.

\begin{lemma}\label{lemma:hollow_simplices_concentration}
    Fix $i \in [2, d]$ with $d \leq \log_2\ln n + \log_2\log_2\ln n - \log_2(-\ln p^{*})$ and $a \geq 0$, then:
    \begin{align*}
    \PP(H_i = 0) \leq  \exp \left( - n^{\log_2 \log_2 \ln n(1 + o(1))}\right)  
    &&
    \PP(|H_i - \mu_i| \geq  a) \leq 2\exp\left(-\frac{2a^2}{(i+1)^{i+2}\binom{n}{i-1}\binom{n}{i+1}}\right).
\end{align*}

In particular, $\lim_{n \to \infty} \PP(I(\multixnp) = d) = 1$ as long as $p^{*}$ is a positive constant.
\end{lemma}

\begin{proof}
    Fix $i \in [2, d]$ and $\epsilon > 0$. For each $s \in \binom{[n]}{i+1}$ we define $X_s = \ind{s \text{ spans a hollow } i\text{-simplex}}$. Then 
\begin{align*}
    H_i = \sum_{s \in \binom{[n]}{i+1}}X_s && \text{and hence} && \mu_i \coloneqq \EE(H_i) = \binom{n}{i+1} P_i.
\end{align*}

{
Write $\Delta_i = \sum_{s \neq t \in \binom{[n]}{i+1}} \ind{|s \cap t| \geq 2} \EEbracket{X_s X_t}$. As long as $\Delta_i \geq \mu_i$, we can bound the denominator in Janson's inequality \cite[Theorem 1]{janson1990poisson} by $\mu_i + \Delta_i \leq 2 \Delta_i$. Hence we get
\begin{align*}
    \PP(H_i = 0) \leq \exp \left(-\frac{\mu_i^2}{4 \Delta_i} \right) &\leq \exp \left( - \frac{\binom{n}{i+1}^2 P_i^2}{4\sum_{m = 2}^{i} \binom{n}{i+1} \binom{i+1}{m} \binom{n}{i+1-m} P_i^2 \prod_{j=1}^{i-1} p_j^{-\binom{m}{j+1}}} \right) \\
    &= \exp \left( - \frac{\binom{n}{i+1}}{4\sum_{m = 2}^{i} \binom{i+1}{m} \binom{n}{i+1-m} \prod_{j=1}^{i-1} p_j^{-\binom{m}{j+1}}} \right).
\end{align*}
}

{
We are interested in the case when $i = \log_2\ln n + \log_2\log_2\ln n +c$ for some constant $c \in \RR$. In this case, by examining the ratio of two consecutive summands in the denominator sum, we can see that the largest summand will correspond to $m = i$. Therefore, we get the bound
\begin{align*}
     \PP(H_i = 0) \leq \exp \left( - \frac{\binom{n}{i+1}(p^*)^{2^i - i - 1}}{4(i+1)^2 n } \right).
\end{align*}

\gr{A}s long as $c \leq -\log_2(-\ln p^{*})$, the bound is at most $\exp \left( - n^{\log_2 \log_2 \ln n(1 + o(1))}\right)$. The first inequality follows by recalling that $H_k = 0$ implies $H_j = 0$ for all $j \geq k$.
}

{
For concentration, we apply \gr{the} large deviation result \gr{Corollary 2.2 in} \cite{janson2004large}. Fix a particular $s \in \binom{[n]}{i+1}$. Note that $X_s$ and $X_t$ are dependent iff $|s \cap t| \geq 2$. So the number of variables in the sum of $H_i$ that $X_s$ depends on is $\sum_{m=2}^{i+1} \binom{i+1}{m}\binom{n-i-1}{i+1-m} \leq (i+1)^{i+1}(i-1)\binom{n}{i-1} \gr{\leq}   \, (i+1)^{i+2}\binom{n}{i-1}$. Now applying the large deviation result \gr{gives} the second inequality.
}
\end{proof}

\begin{theorem}\label{theorem:mle}
Let $i \in [2, d]$ with $d \leq \log_2\ln n + \log_2\log_2\ln n - \log_2(-\ln p^{*})$ and $\epsilon > 0$. Then
\begin{enumerate}
    \item \[p_i\left(1 - \exp \left( - n^{\log_2 \log_2 \ln n(1 + o(1))}\right)   \right)\leq \EE(\hat{p}_i) \leq p_i;\]
    \item \begin{align*}
        \Var(\hat{p}_i) &\leq 2p_i(1-p_i) \left( (i+1)^{i+1}n^{-1 + \epsilon} (p^{*})^{-2^{i+1} + i + 3} + \exp\left(-\frac{2n^{2\epsilon}(i-1)^{i-1}}{(i+1)e^{2i}}\right) \right)\\ &+ p_i^2 \exp \left( - n^{\log_2 \log_2 \ln n(1 + o(1))}\right).
    \end{align*}
\end{enumerate}
\end{theorem}

In particular, for all $i \in [2,d]$, $\lim_{n \to \infty} \EE(\hat{p}_i) = p_i$,  and for $i < \log_2 \ln n + \log_2(1-\epsilon) - \log_2(-\ln(p^*))$ we have $\lim_{n \to \infty} \Var(\hat{p}_i) = 0$ as long as $p^{*}$ is a positive constant.
}

\begin{proof}
To prove the first part of the theorem we use the law of total expectation and fact that $S_i$ is conditionally binomial, \gr{to obtain}
\begin{align*}
    \EE(\hat{p}_i) = \sum_{h_i = 1}^{\binom{n}{i+1}} \frac{h_ip_i}{h_i} \PP(H_i = h_i) = p_i \PP(H_i \neq 0) .
\end{align*}
Applying the bounds from Lemma \ref{lemma:hollow_simplices_concentration} the result follows.
\gr{Turning to} the variance of the estimator,
\begin{align*}
    \Var(\hat{p}_i) &= \EE(\hat{p}_i^2) - \EE(\hat{p}_i)^2 =\sum_{h_i = 1}^{\binom{n}{i+1}} \left( \frac{p_i(1-p_i)}{h_i} + p_i^2\right) \PP(H_i = h_i) - p_i^2 \PP(H_i \neq 0)^2\\
    &=\sum_{h_i = 1}^{\binom{n}{i+1}} \frac{p_i(1-p_i)}{h_i} \PP(H_i = h_i) + p_i^2\PP(H_i \neq 0) (1 - \PP(H_i \neq 0) ).
\end{align*}

{For the second term we can immediately use the first inequality in Lemma \ref{lemma:hollow_simplices_concentration}. For the first term, we use the large deviation inequality from Lemma \ref{lemma:hollow_simplices_concentration}. Fix \tat{$\epsilon \in (0,1)$} and write $a = n^{i + \epsilon}$. Then we have}
{
\begin{align*}
    &\sum_{h_i = 1}^{\binom{n}{i+1}} \frac{p_i(1-p_i)}{h_i} \PP(H_i = h_i) \\
    &= \sum_{h_i = 1}^{\mu_i - a} \frac{p_i(1-p_i)}{h_i} \PP(H_i = h_i) + \sum_{h_i = \mu_i - a + 1}^{\mu_i + a - 1} \frac{p_i(1-p_i)}{h_i} \PP(H_i = h_i) + \sum_{h_i = \mu_i + a}^{\binom{n}{i+1}} \frac{p_i(1-p_i)}{h_i} \PP(H_i = h_i)\\
    &\leq p_i(1-p_i)\left(\sum_{h_i = 1}^{\mu_i - a} \PP(H_i = h_i) + \sum_{h_i = \mu_i + a}^{\binom{n}{i+1}} \PP(H_i = h_i)\right) + \sum_{h_i = \mu_i - a + 1}^{\mu_i + a - 1} \frac{p_i(1-p_i)}{\mu_i - a} \\
    &\leq 2p_i(1-p_i) \exp\left(-\frac{2n^{2\epsilon}(i-1)^{i-1}}{(i+1)e^{2i}}\right) + 2n^{i + \epsilon} \frac{p_i(1-p_i)}{\mu_i - \frac{1}{2}\mu_i}.
\end{align*}
To finish the proof note that
$\mu_i^{-1} n^{i + \epsilon} \leq  (i+1)^{i+1}n^{-1 + \epsilon} (p^{*})^{-2^{i+1} + i + 3} .$
}
\end{proof}

\begin{remark}
Note that if $i=1$ then $\hat{p}_i$ is just the MLE of the parameter $p$ in a binomial distribution, which is very well studied and known to be unbiased and consistent. Also, \cite[Theorem 3.1]{farber2020random} shows that there is some constant $c$ such that for $i \geq \log _2 \ln n+\log _2 \log _2 \ln n - c$ we have no $i$-simplices with high probability as long as the parameter $\mathbf{p} \in (0,1)^{n-1}$ does not depend on $n$. In light of this fact we can see that in Theorem \ref{theorem:mle} the upper bound on $i$ is optimal up to an additive constant {in the case of the expectation. In the case of the variance, \gr{the bound is} off by a factor of $\log_2 \log_2 \ln n$, which is likely because the large deviation inequality is sub-optimal for large $i$.}
\end{remark}

\begin{lemma}\label{lemma:mle_covariance}
{
    For any $2 \leq i < j \leq d \leq \log_2\ln n + \log_2\log_2\ln n - \log_2(-\ln p^{*})$ we have
    \begin{align*}
        -2p_ip_j\binom{n}{i+1}^{\frac{1}{2}}\binom{n}{j+1}^{\frac{1}{2}}&\exp \left( - n^{\log_2 \log_2 \ln n(1 + o(1))}\right)   \leq \\ 
        &\Cov(W_i, W_j) \leq p_ip_j\binom{n}{i+1}^{\frac{1}{2}}\binom{n}{j+1}^{\frac{1}{2}} \exp \left( - n^{\log_2 \log_2 \ln n(1 + o(1))}\right)  .
    \end{align*}
    Also, for any $j \leq d$ we have $\Cov(W_1, W_j) = 0$. Moreover, we have $\lim_{n \to \infty}\Cov(W_i, W_j) = 0$ as long as $p^{*}$ is a positive constant.
   }
\end{lemma}

\begin{proof}
We will use a shorthand $F_n = \binom{n}{i+1}^{\frac{1}{2}}\binom{n}{j+1}^{\frac{1}{2}}$. Also note that {for $i \neq j$,} $S_i | H_i = x$ and $S_j | H_j = y$ are independent Binomial variables. Hence
    \begin{align*}
        \Cov(W_i, W_j) &= F_n \Cov(\hat{p}_j, \hat{p}_j) \\
        &= \sum_{x = 1}^{\binom{n}{i+1}} \sum_{y = 1}^{\binom{n}{j+1}} \EEbracket{\frac{S_iS_j}{xy}}\PP(H_i = x \cap H_j = y) - p_ip_j \PP(H_i \neq 0) \PP(H_j \neq 0)\\
        &=F_n\sum_{x = 1}^{\binom{n}{i+1}} \sum_{y = 1}^{\binom{n}{j+1}} \frac{p_ixp_jy}{xy}\PP(H_i = x \cap H_j = y) - F_np_ip_j \PP(H_i \neq 0) \PP(H_j \neq 0)\\
        &= F_n p_ip_j \left(\PP(H_i \neq 0 \cap H_j \neq 0) - \PP(H_i \neq 0) \PP(H_j \neq 0) \right).
    \end{align*}

    Note that if $i = 1$, then $H_i$ is the number of vertices so $\PP(H_i \neq 0) = 1$ and the covariance is exactly zero.

{    Now from Lemma \ref{lemma:hollow_simplices_concentration}, 
\[1 - \exp \left( - n^{\log_2 \log_2 \ln n(1 + o(1))}\right)   \leq \PP(H_i \neq 0) \leq 1.\]}

 Using this bound together with  $\PP(A) + \PP(B) - 1 \gr{\leq} \PP(A \cap B) \leq \min(\PP(A), \PP(B))$ gives the result.
\end{proof}

\begin{theorem}\label{theorem:mle_normality}
 Let $i \in [2, d]$ with $d \leq \log_2\ln n + \log_2\log_2\ln n - \log_2(-\ln p^{*})$. Also, let \tat{$W_i = \sqrt{\binom{n}{i+1}P_i}(\hat{p}_i - p_i)$ and $Z_i \sim \normal{0}{p_i(1-p_i)}$}. Then
   \begin{align*}
       \dks{W_i}{Z_i} &\leq \frac{\sqrt{10}+3}{3 \sqrt{\pi}}(p_i(1-p_i)P_i)^{-\frac{1}{2}}\binom{n}{i+1}^{-\frac{1}{2}} + \tat{ 2(n-i)^{-2}P_i^{-1} (i+1)^{\frac{i+6}{2}} \ln^{\frac{1}{2}}\binom{n}{i+1} } + 2\binom{n}{i+1}^{-2} .
   \end{align*}
\end{theorem}
 
  \begin{remark}
      \gr{Theorem \ref{theorem:mle_normality} gives in}
  particular, for $i < \log_2 \ln n + 1 - \log_2(-\ln(p^*))$ we have $\lim_{n \to \infty} \dks{W_i}{Z_i} = 0$ as long as $p^{*}$ is a positive constant.
  This assumption can be relaxed; it suffices that $P_i$ and $p_j, j \le d$ depend on $n$ in such a way that the expression in Theorem \ref{theorem:mle_normality} tends to 0 with $n \rightarrow \infty$ while ensuring that $p_i(1-p_i)$ is bounded away from 0.  
  \end{remark}
  
\begin{proof}
Let $Z \sim \normal{0}{1}$. Recall that $H_i$ is the number of hollow $i$-simplices in $\multixnp$. Define the interval $I_i \coloneqq (\EEbracket{H_i} - a, \EEbracket{H_i} + a)$, where $a>0$ is to be chosen later. Let $h$ be an indicator function of an interval $(-\infty, b)$ for some $b \in \RR$. Set \tat{$c(x,i) = x^{\frac{1}{2}}\binom{n}{i+1}^{-\frac{1}{2}}P_i^{-\frac{1}{2}}(p_i(1-p_i))^{-\frac{1}{2}}$}. Recall that given two real random variables $X, Y$ we write $\dks{X}{Y} \coloneqq \sup_{x \in \RR} \abs{\PP(X \leq x) - \PP(Y \leq x)}$ for the Kolmogorov-Smirnov distance.

To prove the theorem we condition on $H_i$. We use large deviation bounds Lemma \ref{lemma:hollow_simplices_concentration} and inside the high probability region use a classic CLT comparing a binomial to a normal. For comparing a binomial variable to a normal we use \cite[Theorem 1]{schulz2016optimal} and for comparing a normal to a normal we use \cite[Section 6.2]{ley2017steins}.

\begin{align*}
    &\EEbracket{\abs{h(W_i) - h(Z_i)}} \\
    &= \sum_{x \in I_i}\EEbracket{\abs{h(W_i) - h(Z_i)} | H_i = x}\PP(H_i = x) + \EEbracket{\abs{h(W_i) - h(Z_i)} | H_i \notin I_i}\PP(H_i \notin I_i) \\
    &\leq \sum_{x \in I_i}\dks{W_i | H_i = x}{Z_i} \PP(H_i = x) + 2\exp\left(-\frac{2a^2}{(i+1)^{i+2}\binom{n}{i-1}\binom{n}{i+1}}\right)\\
    &\leq \sum_{x \in I_i}\dks{c(x,i)W_i | H_i = x}{c(x,i)Z_i} \PP(H_i = x) + 2\exp\left(-\frac{2a^2}{(i+1)^{i+2}\binom{n}{i-1}\binom{n}{i+1}}\right)\\
    &\leq \sum_{x \in I_i} \left(\dks{c(x,i)W_i | H_i = x}{Z} + \dks{c(x,i)Z_i}{Z} \right)  \PP(H_i = x) + 2\exp\left(-\frac{2a^2}{(i+1)^{i+2}\binom{n}{i-1}\binom{n}{i+1}}\right) 
\end{align*}
Note that $c(x,i)W_i | H_i = x$ is a binomial random variable with $n = x$ and $p = p_i$, which has been centered and scaled to have unit variance. Therefore, we can apply a classic CLT to bound the Kolmogorov-Smirnov distance between a standard normal and a centered and scaled binomial \cite[Theorem 1]{schulz2016optimal}. This gives
the bound
\begin{align*}
    \sum_{x \in I_i} \dks{c(x,i)W_i | H_i = x}{Z}\PP(H_i = x) &\leq \sum_{x \in I_i} \frac{\sqrt{10}+3}{6 \sqrt{2 \pi}} \cdot \frac{p_i^2+(1-p_i)^2}{\sqrt{x p_i (1-p_i)}} \PP(H_i = x) \\
    &\leq \frac{\sqrt{10}+3}{6 \sqrt{2 \pi}} \cdot \frac{p_i^2+(1-p_i)^2}{\sqrt{(\EEbracket{H_i} - a) p_i (1-p_i)}} \PP(H_i \in I_i).
\end{align*}

Now to \gr{compare} two normal variables, namely, $c(x,i)Z_i$ and $Z$, we note that \tat{$\Var(c(x,i)Z_i) = \frac{x}{\EEbracket{H_i}}$} \gr{and $\Var(Z)=1$,} and we use  \cite[Section 6.2]{ley2017steins} \gr{(with $f_1=f_2=1$) to obtain}
\tat{\begin{align*}
    \sum_{\substack{x \in I_i }}\dks{c(x,i)Z_i}{Z} \PP(H_i = x) 
    &\leq \sum_{x \in I_i}    \left|\frac{x}{\EEbracket{H_i}}-1\right| \sqrt{\frac{\EEbracket{H_i}}{x}} \PP(H_i = x) \\
     &\leq  \gr{ \max \left\{ \left|\frac{\EEbracket{H_i} + a}{\EEbracket{H_i}}-1\right|, \left|\frac{\EEbracket{H_i} -  a}{\EEbracket{H_i}}-1\right|\right\} }  \sqrt{\frac{\EEbracket{H_i}}{\EEbracket{H_i} - a}} \PP(H_i \in I_i ) \\
     &\leq \frac{a}{\EEbracket{H_i} - a}.
\end{align*}}
Putting everything together, we get
\begin{align*}
    \EEbracket{\abs{h(W_i) - h(Z_i)}} 
     \leq \frac{\sqrt{10}+3}{6 \sqrt{\pi}} \frac{p_i^2+(1-p_i)^2}{\sqrt{\EEbracket{H_i} p_i (1-p_i)}}  + \tat{\frac{a}{\EEbracket{H_i} - a}} + 2\exp\left(-\frac{2a^2}{(i+1)^{i+2}\binom{n}{i-1}\binom{n}{i+1}}\right) .
\end{align*}

Now \gr{for a bound in Kolmogorov distance w}e take \gr{the} supremum over \gr{all half-open interval indicator functions} $h$ and recall that $\EEbracket{H_i} = \binom{n}{i+1}P_i$. \tat{Picking $a = \left( \binom{n}{i+1}\binom{n}{i-1} (i+1)^{i+2} \ln \binom{n}{i+1} \right)^{\frac{1}{2}}$ \gr{to approximately balance the last two summands gives the assertion.}}
\end{proof}

\subsection{\gr{A} goodness of fit test}\label{subsection:testing}

Theorem  \ref{theorem:crit_lexi_approx} provides a bound on the normal approximation error of critical simplex counts with respect to the lexicographical matching. In network analysis goodness of fit tests based on asymptotic distributions of random graph statistics are \gr{available, see for example} \cite{lei2016goodness,ouadah2020degree}. Inspired by this approach, we develop \gr{a} goodness of fit test based on the asymptotic distribution of the critical simplex counts. It is known that in certain geometric random simplicial complexes the lexicographical matching behaves differently compared to the $\multixnp$ model, see \cite[proof of Theorem 5.1]{kahle2011random}. With this in mind, we investigate if the statistic can distinguish between geometric models and $\multixnp$. 

The starting point for the testing procedure is an observed simplicial complex $\cplx$. The null hypothesis is that the observed data is a sample from $\multixnp[n][\hat{\mathbf{p}}_0]$, where $\hat{\mathbf{p}}_0$ is the observed value of the MLE of $\mathbf{p}$ {as given in Equation \ref{equation:mle}}. 
\gr{For $T = (T_1, \ldots, T_n)$ the vector of critical simplex counts in $\multixnp$, we calculate the standardised statistics 
$W_i = W(T_i) = {\Var(T_i)}^{-\frac{1}{2}}(T_i - {\EEbracket{T_i}})$ for  $i \in [n]$. Theorem \ref{theorem:main_CLT} combined with Lemma  \ref{lemma:mle_covariance} give that $W = (W_1, \ldots, W_n)  \approx \mvn{0}{\Sigma}$, where $\Sigma$ is the limiting, diagonal covariance matrix of $W$.  Hence, in the $\multixnp$ model, $W^T \Sigma^{-1} W $ is approximately chi-square distributed with $n$ degrees of freedom. Based on this asymptotic result, as} a proof of concept we propose and empirically study the following simple testing procedure that is described in Algorithm \ref{algorithm:test}:
\begin{enumerate}
    \item We observe a simplicial complex $\cplx$ and compute the corresponding value $\hat{\mathbf{p}}_0$ of the estimator $\hat{\mathbf{p}}$ {from Equation \ref{equation:mle}}.
    \item We compute the observed value $\mathbf{t} = (t_1, t_2, \ldots, t_n)$ of \gr{the critical simplex counts}
    $T = (T_1, T_2, \ldots, T_n)$. 
    \item Based on $\hat{\mathbf{p}}_0$ we estimate $\EEbracket{T_i}$ and $\Var(T_i)$ by replacing $\mathbf{p}$ by $\hat{\mathbf{p}}_0$ in the explicit formulas. We denote these estimators $\widehat{\EEbracket{T_i}}$ and $\widehat{\Var(T_i)}$. 
    \item We \gr{calculate} $\mathbf{w} = (w_1, w_2, \ldots, w_n)$ where $w_i = \widehat{\Var(T_i)}^{-\frac{1}{2}}(t_i - \widehat{\EEbracket{T_i}})$ for all $i \in [n]$.
    \item We reject the null at \tat{approximate} significance {level} $\alpha \in (0,1)$ iff $w^T\Sigma^{-1}w > \chi^2_{n}(1 - \alpha)$, where $\chi^2_{n}(p)$ is the quantile function of the chi-squared distribution with $n$ degrees of freedom.
\end{enumerate}

\begin{algorithm}[H]
	\caption{Goodness-of-fit test}\label{algorithm:test}
	\begin{algorithmic}[1]
		\Procedure {PerformTest}{$\cplx$, $\alpha$, $T$, $\Sigma$} 
		\Comment{$\cplx$ is a finite observed simplicial complex, $\alpha \in (0,1)$ is significance {level}, $T$ is a test statistic, $\Sigma$ is the limiting  covariance of a centered and normalised version of $T$}
		\State $\hat{\mathbf{p}}_0 \coloneqq \Call{ComputeMle}{\cplx}$
		\Comment{Computes the observed value $\hat{\mathbf{p}}_0$ of the MLE $\hat{\mathbf{p}}$}
            \State $(t_1, t_2, \ldots, t_k) \coloneqq \Call{ComputeStatistic}{\cplx, T}$
		\Comment{Computes the observed value of a chosen $T$}
            \ForAll {$i \in [k]$}
            \State $\widehat{\EEbracket{T_i}} \coloneqq \Call{EstimateMean}{\hat{\mathbf{p}}_0}$ \Comment{Estimates the mean by replacing $\mathbf{p}$ with $\hat{\mathbf{p}}_0$}
            \State $\widehat{\Var(T_i)} \coloneqq \Call{EstimateVariance}{\hat{\mathbf{p}}_0}$ \Comment{Estimates the variance by replacing $\mathbf{p}$ with $\hat{\mathbf{p}}_0$}
            \State $w_i \coloneqq \widehat{\Var(T_i)}^{-\frac{1}{2}}(t_i - \widehat{\EEbracket{T_i}})$
            \Comment{Normalises the statistic}
            \EndFor
            \State $w \coloneqq (w_1, w_2, \ldots, w_k)$
            \If {$w^T\Sigma^{-1}w > \chi^2_{k}(1 - \alpha)$} \Comment{$\chi^2_{k}(p)$ is the quantile function of the $\chi^2_{k}$ distribution}
            \State \Return False \Comment{The null is rejected}
            \Else
            \State \Return True \Comment{Insufficient evidence to reject the null}
            \EndIf
		\EndProcedure
		\Statex

	\end{algorithmic}
\end{algorithm}

\begin{remark}
This testing procedure can be improved in the following ways. \gr{First, it ignores the parameter estimation error, second, it ignores the estimation error of the mean and variances, and third, it ignores the normal approximation error.} 
By using Theorem \ref{theorem:mle_normality}, one could account for the parameter estimation error. It is possible to account for estimation error when estimating the mean and the variance of the test statistic as well. Moreover, Theorem \ref{theorem:crit_lexi_approx} gives explicit bounds on \tat{the quality of} the normal approximation and hence it is also possible to account for this error. These improvements are outside the scope of this paper and are deferred to future work.
\end{remark}

\subsubsection{Empirical simulation study}

Here we empirically study the goodness of fit test to distinguish between $\multixnp$ and a random geometric simplicial complex model, where geometry is present only in a particular dimension. We start by defining a soft geometric random simplicial complex model.

\begin{definition}
    Let $\mathcal{X}_n$ be a set of $n$ points that are i.i.d. in the unit $d$-cube $[0,1]^d \subseteq \RR^d$. For $i \in [n-1]$ let $\phi_i: ([0,1]^d)^{i + 1} \to [0,1]$ be a symmetric function. We inductively define a \textbf{soft random geometric simplicial complex} $\mathbb{X}(\mathcal{X}_n, \phi_1, \phi_2, \ldots \phi_{n-1})$ as the random simplicial complex with vertex set $\mathcal{X}_n$ {such that} any $i$-simplex $s = \set{s_1, s_2, \ldots, s_{i+1}} \subseteq \mathcal{X}_n$ is included with probability $\phi_i(s_1, s_2, \ldots, s_{i+1})$ {if} every proper subset of $s$ is already in the complex. Note that if we take the functions $\phi_i$ to be constant functions that do not depend on the vertex locations for all $i \in [n-1]$, then we recover the $\multixnp$ model.
\end{definition}

In the simulation study we compare $\multixnp$ with three geometric models by defining an interpolation process that starts with $\multixnp$ and finishes at each of the models. We are interested to see at which point in the interpolation there is enough geometry in the model for the statistical test to pick it up. We use the following three models:
\begin{enumerate}
    \item Let $\mathcal{X}_{150}$ be a set of $150$ points that are i.i.d. in the unit $7$-cube $[0,1]^7 \subseteq \RR^7$. Given four points $x_1, x_2, x_3, x_4 \in [0,1]^7$ define $A(x_1,x_2,x_3,x_4)$ to be the {volume} of the tetrahedron defined by the four points. The model is $\mathbb{X}(\mathcal{X}_{150}, \phi_1, \phi_2, \ldots \phi_{149})$ where $\phi_1(x_1,x_2) = 0.5$, $\phi_2(x_1,x_2,x_3) = 0.5$, $\phi_3(x_1,x_2,x_3,x_4) = \ind{A(x_1,x_2,x_3,x_4) \leq 0.09}$ and $\phi_i(x_1, x_2, \ldots, x_{i+1}) = 0$ for all $i \neq 1,2,3$. This is just a $3$-dimensional $\multixnp$ model with the exception that only tetrahedra of sufficiently small volume appear. {We later refer to this model as the \textbf{tetrahedron model}.}
    \item Let $\mathcal{X}_{75}$ be a set of $75$ points that are i.i.d. in the unit $3$-cube $[0,1]^3 \subseteq \RR^3$. Given three points $x_1, x_2, x_3 \in [0,1]^3$ define $A(x_1,x_2,x_3)$ to be the area of the triangle defined by the three points. The model is $\mathbb{X}(\mathcal{X}_{75}, \phi_1, \phi_2, \ldots \phi_{74})$ where $\phi_1(x_1,x_2) = 0.5$, $\phi_2(x_1,x_2,x_3) = \ind{A(x_1,x_2,x_3) \leq 0.09}$ as well as $\phi_i(x_1, x_2, \ldots, x_{i+1}) = 0$ for all $i \neq 1,2$. This is just a $2$-dimensional $\multixnp$ model with the exception that only triangles of sufficiently small area appear. {We later refer to this model as the \textbf{triangle model}.}
    \item Let $\mathcal{X}_{75}$ be a set of $75$ points that are i.i.d. in the unit $3$-cube $[0,1]^3 \subseteq \RR^3$. The model here is $\mathbb{X}(\mathcal{X}_{75}, \phi_1, \phi_2, \ldots \phi_{74})$ where $\phi_1(x_1,x_2) = \ind{\norm{x_1 - x_2}_2 \leq 0.4924}$ and $\phi_i(x_1, x_2, \ldots, x_{i+1}) = 0.5$ for all $i \neq 1$. {The distance threshold is chosen so that the edge density of the geometric graph is $0.5$.} This model is just a classical geometric random graph with higher simplices filled in combinatorially just like in $\multixnp$. {We later refer to this model as the \textbf{edge model}.}
\end{enumerate}

In the case of the \tat{tetrahedron} model, we \tat{define a sequence of models interpolating} between the geometric model and $\multixnp$ \tat{as follows}. Given $0 \leq \epsilon_1 \leq \epsilon_2$ we define the model $\mathbb{X}(\mathcal{X}_{150}, \phi_1, \phi_2, \phi_3^{(\epsilon_1, \epsilon_2)}, \ldots, \phi_{149})$ where we set $\phi_1(x_1,x_2) = 0.5$, $\phi_2(x_1,x_2,x_3) = 0.5$, $\phi_3(x_1,x_2,x_3,x_4) = \ind{A(x_1,x_2,x_3,x_4) \leq \epsilon_1} + 0.5 \times \ind{A(x_1,x_2,x_3,x_4) \in (\epsilon_1, \epsilon_2)}$ and $\phi_i(x_1, x_2, \ldots, x_{i+1}) = 0$ for all $i \neq 1,2,3$. That is, we connect any two vertices independently with probability $0.5$, and any triangle becomes a $2$-simplex, also independently, with probability $0.5$. Finally, any hollow $3$-simplex is filled in with probability 1 if its volume is at most $\epsilon_1$ and, independently, with probability $0.5$ if its volume is between $\epsilon_1$ and $\epsilon_2$. No hollow $3$-simplex with volume at least $\epsilon_2$ gets filled in. Note that if we have $\epsilon_1 = \epsilon_2 = 0.09$, then we recover the geometric model and if $\epsilon_1 = 0$, $\epsilon_2 = 1$, then we recover the $\multixnpnotbold[150][(0.5, 0.5, 0.5, 0)]$ model. We start the interpolation with $\epsilon_1 = 0$ and $\epsilon_2 = 0.5$, which is very close to the $\multixnpnotbold[150][(0.5, 0.5, 0.5, 0)]$ model. We increase $\epsilon_1$ and decrease $\epsilon_2$ at each step until we reach $\epsilon_1 = \epsilon_2 = 0.09$. We sample 100 simplicial complexes at each step of the interpolation and report the number of times the null hypothesis is not rejected at $0.95$ significance {level}. 

In the \tat{edge} model, we define a very similar interpolation to the \tat{tetrahedron model interpolation}. Given $0 \leq \epsilon_1 \leq \epsilon_2$ we define $\mathbb{X}(\mathcal{X}_{75}, \phi_1, \phi^{(\epsilon_1, \epsilon_2)}_2, \ldots \phi_{74})$ where $\phi_1(x_1,x_2) = 0.5$, $\phi_2(x_1,x_2,x_3) = \ind{A(x_1,x_2,x_3) \leq \epsilon_1} + 0.5 \times \ind{A(x_1,x_2,x_3) \in (\epsilon_1, \epsilon_2)}$ and $\phi_i(x_1, x_2, \ldots, x_{i+1}) = 0$ for all $i \neq 1,2$. We start the interpolation at $\epsilon_1 = 0$ and $\epsilon_2 = 0.66$, which is close to the $\multixnpnotbold[75][(0.5, 0.5, 0)]$ model. We increase $\epsilon_1$ and decrease $\epsilon_2$ at each step until we reach $\epsilon_1 = \epsilon_2 = 0.09$. We sample 100 simplicial complexes at each step of the interpolation and report the number of times the null hypothesis is not rejected at $0.95$ \tat{approximate} significance {level}.

Finally, in the \tat{edge} model we perform an analogous interpolation again. Given $0 \leq \epsilon_1 \leq \epsilon_2$ we define $\mathbb{X}(\mathcal{X}_{75}, \phi_1, \phi^{(\epsilon_1, \epsilon_2)}_2, \ldots \phi_{74})$ where $\phi_1(x_1,x_2) =\ind{\norm{x_1 - x_2}_2 \leq \epsilon_1} + 0.5 \times \ind{\norm{x_1 - x_2}_2 \in (\epsilon_1, \epsilon_2)}$ as well as $\phi_i(x_1, x_2, \ldots, x_{i+1}) = 0.5$ for all $i \neq 1$. We start the interpolation at $\epsilon_1 = 0$ and $\epsilon_2 = \sqrt{3}$, which is exactly the $\multixnpnotbold[75][(0.5, 0.5, 0.5 \ldots, 0.5)]$ model. We increase $\epsilon_1$ and decrease $\epsilon_2$ at each step until we reach $\epsilon_1 = \epsilon_2 = 0.4924$. 

\gr{For each of these three perturbation models we} sample 100 simplicial complexes at each step of the interpolation and report the number of times the null hypothesis \gr{of an $\bf{X}(n, {\bf{p}})$ model} is not rejected \gr{(the {\it number of passes})} at $0.95$ \tat{approximate} significance {level}.  \gr{Ideally the number of passes should be 0 when  $\Delta \epsilon =0$, corresponding to the purely geometric model,  and increase sharply with increasing $\Delta \epsilon,$ to  100 when $\Delta \epsilon =1$, corresponding to the $\bf{X}(n, {\bf{p}})$ model.} 

\begin{figure}[H]
		\centering
		\includegraphics[width=0.8\textwidth]{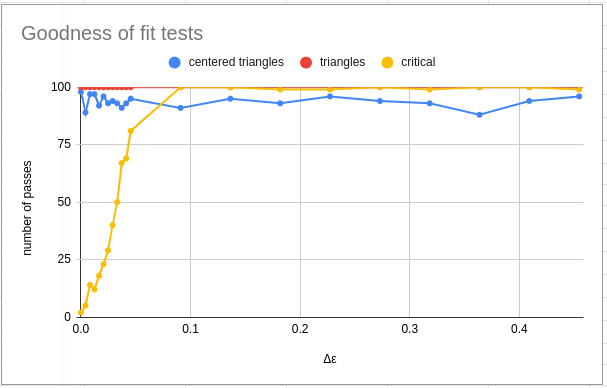}
		\caption{Goodness of fit testing results of the tetrahedron model. The $x$-axis corresponds to $\Delta\epsilon \coloneqq \epsilon_2 - \epsilon_1$.}
		\label{fig:dim3_tests}
\end{figure}

\begin{figure}[H]
		\centering
		\includegraphics[width=0.8\textwidth]{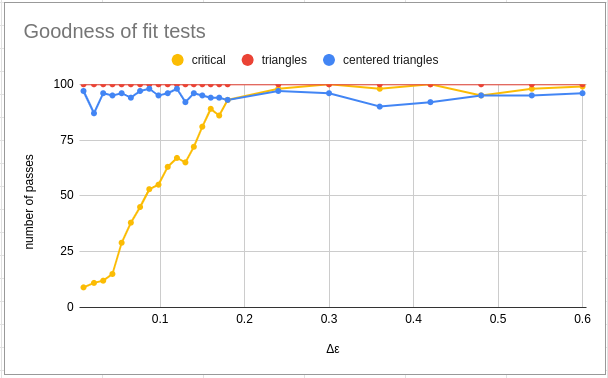}
		\caption{Goodness of fit testing results of the triangle model. The $x$-axis corresponds to $\Delta\epsilon \coloneqq \epsilon_2 - \epsilon_1$.}
		\label{fig:dim2_tests}
\end{figure}

\begin{figure}[h]
		\centering
		\includegraphics[width=0.8\textwidth]{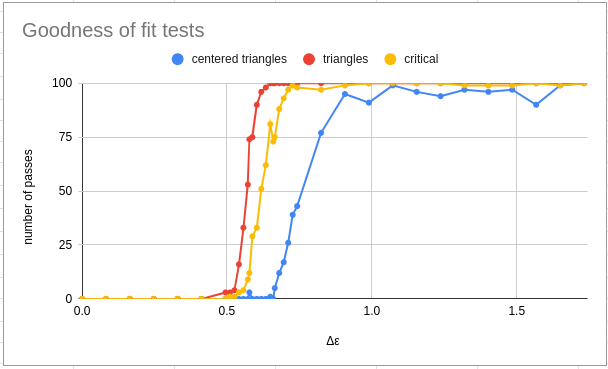}
		\caption{Goodness of fit testing results of the edge model. The $x$-axis corresponds to $\Delta\epsilon \coloneqq \epsilon_2 - \epsilon_1$.}
		\label{fig:dim1_tests}
	\end{figure}

We compare goodness of fit tests based on three statistics: triangle counts, centered triangle counts \cite{kaur2020higher, bubeck2016testing}, and the critical simplex counts.  $\Delta \epsilon$. Because the triangle counts and centered triangle counts only rely on $1$-dimensional information, we expect these statistics to detect the difference only when comparing $\multixnp$ to the last model. \tat{The results for the tetrahedron, triangle, and edge models can be seen in Figure \ref{fig:dim3_tests}, \ref{fig:dim2_tests}, and \ref{fig:dim1_tests} respectively. The data used to generate the respective figures can be found in Table \ref{table:dim3}, \ref{table:dim2}, and \ref{table:dim1} in Appendix \ref{section:appendix_simulation}. }

\paragraph{The simulation results} We \tat{find} that critical simplex counts under the lexicographical matching can differentiate between {the} $\multixnp$ model, which is fully combinatorial, and models which depend on geometric information only in a particular dimension. Even if the model is far from a fully geometric one (like a Vietoris-Rips or \v{C}ech complex of a point process in $\RR^d$ would be), \tat{we can still distinguish it from the $\multixnp$ based on the critical simplex counts}. In dimension 1 it is relatively easy to distinguish between a low-dimensional geometric random graph and $\gnp$ as it has been noted by others \cite{bubeck2016testing, liu2021phase}. In this case, \tat{as seen in Figure \ref{fig:dim1_tests},} the centered triangle counts can tolerate more combinatorial noise and it starts detecting geometry earlier than the number of critical simplices. However, we have noticed {based on the test} simulations that the centered triangle counts are less stable and affected by parameter estimation errors. Quite often we would see less than 90\% passes in the cases where the null hypothesis is true at \tat{approximate} $0.95$ significance {level}. Generalising centered subgraph counts to centered subcomplex counts and investigating their use in goodness-of-fit tests would be a fruitful research project. An advantage of critical simplex counts, beyond their performance as a test statistic, is that when performing persistent homology computations in practice one quite often tries to reduce the computational cost by applying a lexicographical acyclic partial matching and so a statistical test based on the number of critical simplices can be performed at very little added cost.

\section{Probabilistic tools}\label{section:probabilistic_tools}

{In this section we present the probabilistic tools that we use to prove the multivariate CLTs for the subcomplex counts and the number of critical simplices under the lexicographical matching.}

\subsection{Multivariate CLT for generalised $U$-statistics}

\tat{To prove the results in this section we rely on the theorems from \cite{temvcinas2021multivariate}, which are briefly introduced in Appendix \ref{section:appendix_clt} for completeness.}

Assume we have $k$ not necessarily independent collections of independent random elements, where the $i$-th collection is indexed by $\binom{[n]}{i}$, the set of all subsets of $[n]$ of size $i$. That is, the $i$-th collection is given by \[\set{\xi^{(i)}_{\alpha}}_{\alpha \in \binom{[n]}{i}}.\]
Here $i$ is any positive integer not bigger than $n$. For example, one can think of $\xi^{(i)}_{\alpha}$ as an indicator of a hyper-edge $\alpha$ of size $i$ in some random hypergraph model in which the hyperedges occur independently.

Now we consider the sequence of collections \[\set{\xi^{(1)}_{\alpha}}_{\alpha \in \binom{[n]}{1}}, \set{\xi^{(2)}_{\alpha}}_{\alpha \in \binom{[n]}{2}}, \ldots, \set{\xi^{(k)}_{\alpha}}_{\alpha \in \binom{[n]}{k}}.\]

Assume that the random elements $\set{\xi^{(i)}_{\alpha}}_{\alpha \in \binom{[n]}{i}}$ take values in a measurable set $\mathbb{X}^{(i)}$ for all $i \in [k]$ and that there is some measurable set $\mathbb{X}$ such that for all $i \in [k]$ we have $\mathbb{X}^{(i)} \subseteq \mathbb{X}$.

Given a subset $s \subseteq [n]$ of size $m \geq k$ let $\binom{s}{i}$ be the set of $i$-subsets of $s$, which, if needed, we can turn into a sequence by ordering the subsets according to the lexicographical ordering. 

Define a sequence (with respect to the lexicographical ordering on $\binom{s}{i}$) of the random elements \begin{align}\label{equation:xis}
    \mathcal{X}^{(i)}_s = (\xi^{(i)}_{\alpha})_{\alpha \in \binom{s}{i}}.
\end{align}

\begin{definition}
\label{definition:u_statistic_order_k}
Given two integers $1 \leq k, m \leq n$ and a measurable function \[f: \prod_{i=1}^k (\mathbb{X}^{(i)})^{\binom{m}{i}} \to \RR\] define the associated \textbf{generalised $\mathbf{U}$-statistic} of order $k$ to be
\[S_{n,m}^{(k)}(f) = \sum_{s \in \binom{[n]}{m}} f(\mathcal{X}^{(1)}_s, \mathcal{X}^{(2)}_s, \ldots, \mathcal{X}^{(k)}_s).\]
\end{definition}

To clarify this definition, \tat{the aim is} to be able to input all the variables corresponding to subsets of size $i$ of a base subset $s$, which we sum over. There are $\binom{m}{i}$ such subsets and since there is a 1-1 correspondence between the subsets and the random elements, this is how many of those random elements we want to be able to input into our function. Finally we take a product over all $i \in [k]$ because we want the function to be capable of depending of subsets of different size (up to  size $k$).
  
    \begin{example} 
    
    Generalised $U$-statistics of order $k > 2$ covers many random variables of interest in the context of random hypergraphs and random simplicial complexes. Such statistics have been studied in \cite[Section 10.7]{janson1991asymptotic}, \cite{de1996central}, \cite[Chapter 11.3]{janson1997gaussian}.  For example, consider a hypergraph on the vertex set $[n]$ where a hyper-edge $\alpha \subseteq n$ is present independently of all other hyper-edges with probability $p$. Here we assume that  $| \alpha| \ge 2$ so that the number of vertices is fixed. For any $\alpha \subseteq [n]$, such that $|\alpha|> 1$, let $\xi^{(|\alpha|)}_{\alpha}$ be the hyper-edge indicator. Then the sequences \[\set{\xi^{(2)}_{\alpha}}_{\alpha \in \binom{[n]}{2}}, \ldots, \set{\xi^{(k)}_{\alpha}}_{\alpha \in \binom{[n]}{k}}\] are all i.i.d. $\ber{p}$ but they are of different size. Because all vertices a present in this hypergraph model, we set $\xi^{(1)}_{\alpha}=1$ for all $\alpha \in \binom{[n]}{1}$, which makes the first sequence of variables deterministic.  Define a function $f: \RR^3 \times \RR^3 \times \RR \to \RR$ by $$f(x_{1},x_2,x_3,y_{1,2},y_{1,3},y_{2,3},z_{1,2,3}) = z_{1,2,3}(1-y_{1,2})(1-y_{1,3})(1-y_{2,3}).$$ Then the associated $U$-statistic of order 3, $S^{(3)}_{n,3}$, is the number of hyper-edges of size 3, which do not contain any hyper-edge of smaller size in this random hypergraph.
    \end{example}

Fix positive integers $\set{m_i}_{i \in [d]}$ and let $\set{f_i}_{i \in [d]}$ be a collection of measurable functions \[f_i: \prod_{t=1}^k (\mathbb{X}^{(t)})^{\binom{m_i}{t}} \to \RR\] and consider the collection of associated generalised $U$-statistics of order $k$ \[\set{S_{n,m_i}^{(k)}(f_i)}_{i \in [d]}.\]

For any $i \in [d]$ let $\I_i \coloneqq \binom{[n]}{m_i} \times \set{i}$. For $s = (\phi,i) \in \I_i$ define $Y_s = f_i(\mathcal{X}^{(1)}_\phi, \mathcal{X}^{(2)}_\phi, \ldots, \mathcal{X}^{(k)}_\phi)$ and \[X_s = \sigma_i^{-1}\left\{Y_s - \mu_s \right\},\]
where $\mu_s = \EEbracket{Y_s}$ and $\sigma_i = \sqrt{\Var(S_{n,m_i}^{(k)}(f_i))}$, assuming the relevant expectations are finite and variances are all non-zero and finite. 

Let $W_i = \sum_{s \in \I_i} X_s$ and define a random vector $W = (W_1, W_2, \ldots, W_d) \in \RR^d$. We are interested in the distribution of $W$. However, to apply Corollary \ref{corollary:mvn_dissociated_decomp} we need some assumptions.

Let $k'$ be the largest non-negative integer such that the random elements \[\set{\xi^{(1)}_{\alpha}}_{\alpha \in \binom{[n]}{1}}, \set{\xi^{(2)}_{\alpha}}_{\alpha \in \binom{[n]}{2}}, \ldots, \set{\xi^{(k')}_{\alpha}}_{\alpha \in \binom{[n]}{k'}}\] are deterministic. Such $k'$ always exists and if there are no deterministic variables, we will say that $k'=0$. It is convenient to have such $k'$ as a parameter when studying Linial-Meshulam random simplicial complexes; there we take a $(k'-1)$-skeleton that is complete and put in higher simplices with different probabilities. For other classical models like $\xnp$ we usually have $k' = 1$ since there is no randomness at the level of vertices.

\begin{assumption}\label{assumption:ho_ustat}
We use the following two assumptions: 
\begin{enumerate}
    \item For any $i \in [d]$ there is $\alpha_i > 0$ such that for all pairs $s,u \in \I_i$ for which $Y_s, Y_u$ are dependent, we have \[\Cov(Y_s,Y_u) \geq \alpha_i.\]
    \item There is $\beta \geq 0$ such that for any $i,j,l \in [d]$ and any $s \in \I_i, t \in \I_j, u \in \I_l$ we have \[\EE \abs{\left\{ Y_s - \mu_{s} \right\}\left\{ Y_t - \mu_{t} \right\}\left\{ Y_u - \mu_{u} \right\}} \leq \beta\] as well as \[\EE \abs{\left\{ Y_s - \mu_{s} \right\}\left\{ Y_t - \mu_{t} \right\}} \EE \abs{ Y_u - \mu_{u}} \leq \beta.\]
\end{enumerate}
\end{assumption}

The assumptions used here are far from necessary for the normal approximation to hold. However, they hold in the settings studied in this work and adopting them makes the bounds quite convenient. 

\begin{theorem}\label{theorem:ustat_approx}
Let $Z \sim \mvn{0}{\text{\rm Id}_{d \times d}}$ and $\Sigma$ be the covariance matrix of $W$, where $W$ satisfies Assumption \ref{assumption:ho_ustat}.
\begin{enumerate}
    \item Let $h \in \testzero$. Then \[
        \abs{\EE h(W) - \EE h(\Sigma^{\frac{1}{2}}Z)} \leq \abs{h}_3  B_{\ref{theorem:ustat_approx}}  n^{-\frac{k'+1}{2}}.
        \]
        
    \item  \[ \sup _{A \in \mathcal{K}}|\PP(W \in A)-\PP(\Sigma^{\frac{1}{2}}Z \in A)| \leq 2^{\frac{7}{2}} 3^{-\frac{3}{4}}d^{\frac{3}{16}}B_{\ref{theorem:ustat_approx}}^{\frac{1}{4}} n^{-\frac{k'+1}{8}}, \]
\end{enumerate}
where
\[
B_{\ref{theorem:ustat_approx}} = \frac{4}{3}\frac{\beta(k'+1)^{\frac{3}{2}(k'+1)}}{((k'+1)!)^2} (\sqrt{2}m )^{5m + 2 - 3(k'+1)} \sum_{(i,j,k)}^d \frac{1}{\sqrt{\alpha_i\alpha_j\alpha_k}}  
\]
and $m = \max_{i \in [d]} m_i.$
\end{theorem}

\begin{proof}
Fix $i,j \in [d]$ and $s = (\phi,i) \in \I_i$. Then using the definition of $k'$, we see that one choice of dependency neighbourhoods is $\D_j(s) = \setc{(\psi, j) \in \I_j}{|\phi \cap \psi| \geq k' + 1}$. \tat{This follows because} in our setup the sequences of random elements \tat{$\xi_\alpha^{(l)}$} are independent \tat{and so} the \tat{summands} associated to two subsets can only be dependent if they share a random element. By definition of $k'$, the smallest such random element that is not deterministic will be associated to a subset of size $k'+1$. Now we see that Corollary \ref{corollary:mvn_dissociated_decomp} applies.

Looking at the size of $\D_j(s)$, {with \eqref{eq:binom}} we have:
\begin{equation}\label{equation:ho_ustat1}
    |\D_j(s)| = \sum_{m=k'+1}^{\min(m_i,m_j)} \binom{m_i}{m} \binom{n-m_i}{m_j - m} \leq n^{m_j-1-k'} \frac{m_i^{m_i+1}}{(k'+1)!} .
\end{equation}
Let us now calculate a lower bound for the variance using the second part of Assumption \ref{assumption:ho_ustat}, noting that $|\mathbb{I}_i| = \binom{n}{m_i}$:

\begin{align*}
    \sigma_i^2 &= \Var(S_{n,m_i}^{(k)}(f_i)) = \sum_{s \in \I_i} \sum_{u \in D_i(s)} \Cov(Y_s,Y_u)\\
    &= \sum_{s \in \I_i} \sum_{m=k'+1}^{m_i} \sum_{\substack{u \in \binom{[n]}{m_i} \\ |u \cap s| = m}} \Cov(Y_s,Y_u) \geq \sum_{m=k'+1}^{m_i} \binom{n}{m_i} \binom{m_i}{m} \binom{n-m_i}{m_i - m} \alpha_i \\
    &= \sum_{m=k'+1}^{m_i} \binom{n}{2m_i-m} \binom{2m_i - m}{m_i} \binom{m_i}{m}  \alpha_i \geq \binom{n}{2m_i-1-k'} \binom{2m_i-1-k'}{m_i} \binom{m_i}{k'+1}\alpha_i \\
    &\geq \frac{n^{2m_i-1-k'}}{(2m_i-1-k')^{2m_i-1-k'}} \frac{(2m_i-1-k')^{m_i}}{m_i^{m_i}} \frac{m_i^{k'+1}}{(k'+1)^{k'+1}}\alpha_i \\
    &= \frac{n^{2m_i-1-k'}m_i^{k'+1-m_i}}{(2m_i-1-k')^{m_i-1-k'} (k'+1)^{k'+1}}\alpha_i .
\end{align*}

Taking the inequality to the power of $-\frac{1}{2}$ we get:
\[
    \sigma_i^{-1} \leq n^{-m_i+\frac{k'+1}{2}}(2m_i^2-m_i(1+k'))^{\frac{m_i}{2} - \frac{k'+1}{2}}(k'+1)^{\frac{k'+1}{2}}\alpha_i^{-\frac{1}{2}}.
\]

Using the language of Corollary \ref{corollary:mvn_dissociated_decomp} and using the last part of Assumption \ref{assumption:ho_ustat} we can say that $\beta_{ijl} \leq \beta(\sigma_i\sigma_j\sigma_l)$. Moreover we can bound $\alpha_{ij}$ in Corollary \ref{corollary:mvn_dissociated_decomp} by \eqref{equation:ho_ustat1}.
The bound on the quantity $B_{\ref{corollary:mvn_dissociated_decomp}}$ from the corollary is: 
\begin{align*}
    &B_{\ref{corollary:mvn_dissociated_decomp}} = \frac{1}{3} \sum_{(i,j,k)} \abs{\I_i}  \alpha_{ij}  \left(\frac{3\alpha_{ik}}{2} + 2\alpha_{jk}\right)  \beta_{ijk} \\
    &\leq \frac{1}{3} \sum_{(i,j,k)} \frac{n^{m_i}}{m_i!} n^{m_j-1-k'} \frac{m_i^{m_i+1}}{(k'+1)!} \left( \frac{3}{2} n^{m_k-1-k'} \frac{m_i^{m_i+1}}{(k'+1)!} + 2 n^{m_k-1-k'} \frac{m_j^{m_j+1}}{(k'+1)!} \right) \beta_{i,j,k} \\
    &\leq \frac{2\beta(k'+1)^{\frac{3}{2}(k'+1)}}{3((k'+1)!)^2} \sum_{(i,j,k)} n^{m_i+m_j+m_k-2(k'+1)} \frac{m_i^{m_i+1}}{m_i!\sqrt{\alpha_i\alpha_j\alpha_k}} \left( m_i^{m_i+1} + m_j^{m_j+1} \right) n^{-m_i-m_j-m_k+ \frac{3}{2}(k'+1)} M_{i,j,k} \\
    &= \frac{2\beta(k'+1)^{\frac{3}{2}(k'+1)}}{3((k'+1)!)^2} \sum_{(i,j,k)} n^{-\frac{k'+1}{2}}\frac{m_i^{m_i+1}(m_i^{m_i+1} + m_j^{m_j+1})}{m_i!\sqrt{\alpha_i\alpha_j\alpha_k}} M_{i,j,k} \\
    &= \left\{\frac{2}{3}\frac{\beta(k'+1)^{\frac{3}{2}(k'+1)}}{((k'+1)!)^2} \sum_{(i,j,k)} \frac{m_i^{m_i+1}(m_i^{m_i+1} + m_j^{m_j+1})}{m_i!\sqrt{\alpha_i\alpha_j\alpha_k}} M_{i,j,k} \right\} n^{-\frac{k'+1}{2}},
\end{align*}

where $M_{i,j,k} = (2m_i^2-m_i(1+k'))^{\frac{m_i}{2} - \frac{k'+1}{2}}(2m_j^2-m_j(1+k'))^{\frac{m_j}{2} - \frac{k'+1}{2}}(2m_k^2-m_k(1+k'))^{\frac{m_k}{2} - \frac{k'+1}{2}}$.

Now \tat{set} $m = \max_{i \in [d]} m_i$. Then we have $(2m_i^2-m_i(1+k'))^{\frac{m_i}{2} - \frac{k'+1}{2}} \leq (2m^2)^{\frac{m}{2} - \frac{k'+1}{2}}$ and so
\[B_{\ref{corollary:mvn_dissociated_decomp}} \leq \left\{\frac{4}{3}\frac{\beta(k'+1)^{\frac{3}{2}(k'+1)}}{((k'+1)!)^2} (\sqrt{2}m )^{5m + 2 - 3(k'+1)} \sum_{(i,j,k)}^d \frac{1}{\sqrt{\alpha_i\alpha_j\alpha_k}}  \right\} n^{-\frac{k'+1}{2}}.\]
\end{proof}

\section{Subcomplex counts}\label{section:subcomplex_counts}

Let $\cplx$ be a fixed connected simplicial complex on a vertex set $[k+1]$. Let $\mathbf{p} = (p_1, p_2, \ldots, p_{n-1}) \in \RR^{n-1}$ be a vector where for any $i \in [n-1]$ we have $p_i \in [0,1]$. Let $k' \in \NN$ be such that for any $i \in [k']$ we have $p_i = 1$. We are interested in \tat{subcomplex counts of $\cplx$ in $\multixnp$} as defined in Section \ref{section:top_preliminaries}. This is a natural generalisation of subgraph counts.

Let us discuss this variable from a probabilistic perspective. Consider the random hypergraph model from Definition \ref{definition:costa-farber}. For any $k \in [n-1]$ and $\alpha \in \binom{[n]}{k+1}$ let $\xi^{(k+1)}_{\alpha} \sim \ber{p_k}$ be the hyperedge indicators in the random hypergraph model. Hence  the random variables in this collection are independent by construction.
Then the indicator of the $k$-simplex $\alpha$ in the random simplicial complex model can be written as a product of the hyperedge indicators: 
\begin{equation}\label{equation:simplex_ind}
    \ind{\alpha \in \multixnp} = \prod_{i = k'+1}^k \prod_{\substack{\beta \subseteq \alpha \\ |\beta| = i+1}} \xi^{(i+1)}_{\beta}.
\end{equation}
Note that by independence we have \[\EEbracket{\ind{\alpha \in \multixnp}} = \prod_{i=k'+1}^k p_i^{\binom{k+1}{i+1}}.\]

Let $[\cplx]$ be the isomorphism class of $\cplx$ in the set of all simplicial complexes on the vertex set $[k + 1]$. That is, $[\cplx]$ is the set of all simplicial complexes on the vertex set $[k+1]$ that are simplicially isomorphic to $\cplx$. For any $s \in \binom{[n]}{k+1}$,  we write $s = \set{s_1, s_2, \ldots, s_{k+1}}$ such that $s_1 < s_2 < \ldots < s_{k+1}$. Then \tat{we define} $\cplx[s]$ \tat{to be the} complex on the vertex set $s$ such that any $\alpha = \set{s_{i_1}, s_{i_2}, \ldots s_{i_l}} \subseteq s$ is a simplex in $\cplx[s]$ iff $\set{i_1, i_2, \ldots, i_l} \in \cplx$. We are interested in the following variable:
\begin{equation}\label{equation:uncentered_count}
    T_{\cplx} \coloneqq \sum_{s \in \binom{[n]}{k+1}} \sum_{\cplx' \in [\cplx]} \prod_{i = k'+1}^k \prod_{\substack{\alpha \in \cplx'[s] \\ |\alpha| = i+1}} \xi^{(i+1)}_{\alpha}.
\end{equation}

Note that this random variable, \tat{which counts the number of copies of $\cplx$ in $\multixnp$}, can be studied in the framework of generalised $U$-statistics. In particular, if the complex $\cplx$ is of dimension $d$, then $T_\cplx$ is a $U$-statistic of order $d + 1$. 

To fit the random variable \tat{$T_\cplx$} in \tat{the framework of generalised $U$-statistics}, we say that the collection of sequences \[\set{\xi^{(1)}_{\alpha}}_{\alpha \in \binom{[n]}{1}}, \set{\xi^{(2)}_{\alpha}}_{\alpha \in \binom{[n]}{2}}, \ldots, \set{\xi^{(d+1)}_{\alpha}}_{\alpha \in \binom{[n]}{d+1}}\] are the underlying variables from which we can build simplex indicators as it is done in \eqref{equation:simplex_ind}. That is, we have $\xi^{(1)}_{\alpha} = 1$ for any $\alpha \in \binom{[n]}{1}$. For any $1 < k \leq d+1$ and $\alpha \in \binom{[n]}{k+1}$ we have $\xi^{(k+1)}_{\alpha} \sim \ber{p_k}$. All the sequences are independent and within each sequence the variables are i.i.d. Recall that $\mathcal{X}^{(i+1)}_s$ from Equation \ref{equation:xis} is a sequence of the random elements, corresponding to the $(i+1)$-subsets of $s$, ordered lexicographically. We write $\mathcal{X}^{(i+1)}_s[j]$ for the $j$-th element of the sequence. Let $o(\alpha)$ be the position of $\alpha \in \binom{s}{i+1}$ in this lexicographical ordering. Now let us define the associated function $f: \prod_{k=2}^{d+1}(\set{0,1})^{\binom{d+1}{k}} \to \RR$ by
\begin{equation}\label{equation:f_i}
    f(\mathcal{X}^{(1)}_s, \mathcal{X}^{(2)}_s, \ldots, \mathcal{X}^{(d+1)}_s) = \sum_{\cplx' \in [\cplx]} \prod_{i = k'+1}^d \prod_{\substack{\alpha \in \cplx'[s] \\ |\alpha| = i+1}} \mathcal{X}^{(i+1)}_s[o(\alpha)].
\end{equation}

Now it is easy to see that
\[T_{\cplx} = \sum_{s \in \binom{[n]}{k+1}}f(\mathcal{X}^{(1)}_s, \mathcal{X}^{(2)}_s, \ldots, \mathcal{X}^{(d+1)}_s)\] and hence the random variable in question is a \tat{generalised} $U$-statistic of order $d+1$.

\subsection{Moments and the limiting covariance}

Let $e_i(\cplx)$ be the number of $i$-simplices in $\cplx$. Then by independence we have
\[\EEbracket{T_{\cplx}} = \binom{n}{k + 1} |[\cplx]| \prod_{i = 1}^k p_i^{e_i(\cplx)}.\]

For example, if $\cplx = \set{1, 2, 3, 4, 12, 13, 23, 24, 34, 123, 234}$ (two 2-simplices with a common 1-simplex), then $\EEbracket{T_{\cplx}} = 3\binom{n}{4} p_1^{5} p_2^2$. Here 1 and 4 can be exchanged, 2 and 3 can be exchanged, and there are no further simplicial isomorphisms and hence $|[\cplx]| =3$. The following lemma is used in \tat{our} limiting covariance calculations.

\begin{lemma}\label{lemma:binomial_limit}
Let $l,m,a$ be positive integers independent of $n \in \NN$ such that $l, m > \frac{a}{2}$. Also, let $N(n)$ be a sequence of integers such that $N(n) = \Theta(n)$. We write $N$ instead of $N(n)$. Then:
\begin{align*}
    &\lim_{n \to \infty} \binom{N}{l+m-a} \binom{l+m-a}{l} \binom{l}{a}   \left\{\binom{N}{2l-a} \binom{2l-a}{l} \binom{l}{a} \binom{N}{2m-a} \binom{2m-a}{m} \binom{m}{a}   \right\}^{-\frac{1}{2}}  = 1.
\end{align*}
\end{lemma}

\begin{proof}
\begin{align*}
    &\lim_{n \to \infty} \binom{N}{l+m-a} \binom{l+m-a}{l} \binom{l}{a} \left\{\binom{N}{2l-a} \binom{2l-a}{l} \binom{l}{a}\binom{N}{2m-a} \binom{2m-a}{m} \binom{m}{a} \right\}^{-\frac{1}{2}} \\
    &=\lim_{n \to \infty} \frac{N!}{(N - (l+m-a))! } \frac{1}{(l!)^{\frac{1}{2}}(m-a)!} \left(\frac{1}{a!(l-a)!}\right)^{\frac{1}{2}}  \\
    &\left\{\frac{N!}{(N - (2l-a))!} \frac{1}{l!(l-a)!} \frac{N!}{(N - (2m-a))!} \frac{1}{((m-a)!)^2} \frac{1}{a!} \right\}^{-\frac{1}{2}} \\
    &=\lim_{n \to \infty} \frac{\left((N - (2m-a))!(N - (2l-a))!\right)^{\frac{1}{2}}}{(N - (l+m -a))! } =1.
\end{align*}
\end{proof}

Theorem \ref{theorem:subcomplex_counts_limiting_cov} shows that, asymptotically, similarly to the subgraph counts in $G(n,p)$ for a constant $p \in (0,1)$, subcomplex counts are perfectly correlated. 

\begin{theorem}\label{theorem:subcomplex_counts_limiting_cov}
Let $\cplx, \mathcal{M}$ be two simplicial complexes of dimension at least $k' + 1$ on the vertex sets $[l]$ and $[m]$ respectively. Let $\sigma_{\cplx}^2 = \Var(T_{\cplx})$ and $\sigma_{\mathcal{M}}^2 = \Var(T_{\mathcal{M}})$. Then assuming that for all $i$ the parameters $p_i$ stay constant, for any fixed $m,l \in \NN$ we have:
\[\lim_{n \to \infty} \sigma_{\cplx}^{-1} \sigma_{\mathcal{M}}^{-1} \Cov(T_{\cplx}, T_{\mathcal{M}}) = 1.\]
\end{theorem}

\begin{proof}
Assume w.l.o.g. that $m \leq l$.  For $s \in \binom{[n]}{l}$ set the variable \[X_s^{\cplx} = \sum_{\cplx' \in [\cplx]} \prod_{i = k'+1}^k \prod_{\substack{\alpha \in \cplx'[s] \\ |\alpha| = i+1}} \xi^{(i+1)}_{\alpha}, \] so that $T_{\mathcal{L}} = \sum_{s \in \binom{[n]}{l}} X_s^{\cplx} $, and analogously $X_u^{\mathcal{M}}$ for $u \in \binom{[n]}{m}$. Note that $X_s^{\cplx}$  does not depend on $n$. If one would like to relate this expression back to the $U$-statistics framework, one needs to note that if $f$ is a function that gives us $T_\cplx$ as a $U$-statistic, then for any $s \in \binom{[n]}{l}$ we have $X_s^{\cplx} = f(\mathcal{X}^{(1)}_s, \mathcal{X}^{(2)}_s, \ldots, \mathcal{X}^{(l)}_s)$.

Note that for any $s \in \binom{[n]}{l}, u \in \binom{[n]}{m}$ if $|s \cap u| < k' + 2$, then the variables $X_s^{\cplx}$ and $X_u^{\mathcal{M}}$ are independent, since the intersection of the underlying subsets is not large enough to contain a random element of the form $\xi^{(i+1)}_{\alpha}$ for $i + 1 \geq k' + 2$, which would create dependence. Hence we see that \[\Cov(T_{\cplx}, T_{\mathcal{M}}) = \sum_{s \in \binom{[n]}{l}} \sum_{k=k'+2}^{m} \sum_{\substack{u \in \binom{[n]}{m} \\ |s \cap u| = k}} \Cov(X_s^{\cplx}, X_u^{\mathcal{M}}).\] 

The limit of the ratio \tat{is determined by} the ratio of the highest order term from the enumerator and the highest order term from the denominator. Those terms will correspond to $k = k'+2$ because the number of pairs $s \in \binom{[n]}{l}, u \in \binom{[n]}{m}$ where $|s \cap u| = k$ is $\binom{n}{l} \binom{l}{k} \binom{n-l}{m-k} = \binom{n}{l+m-k}\binom{l+m-k}{l} \binom{l}{k}$, which is of the order $n^{l+m-k}$, and the summands $\Cov(X_s^{\cplx}, X_u^{\mathcal{M}})$ do not depend on $n$. We write $L_{k'+1}$ for the number of elements of $[\cplx]$ that contain a particularly chosen $(k'+1)-$simplex and $M_{k'+1}$ for the analogous quantity corresponding to $\mathcal{M}$. 

We use Lemma \ref{lemma:binomial_limit} with $N = n$ and $a = k'+2$ to get:
\begin{align*}
    &\lim_{n \to \infty} \sigma_{\cplx}^{-1} \sigma_{\mathcal{M}}^{-1} \Cov(T_{\cplx}, T_{\mathcal{M}}) = \\
    &\lim_{n \to \infty} \binom{n}{l+m-2-k'} \binom{l+m-2-k'}{l} \binom{l}{k'+2} L_{k'+1}M_{k'+1} \left( p_{k'+1}^{-1} - 1 \right) \prod_{i = k'+1}^{l-1} p_i^{e_i(\cplx) + e_i(\mathcal{M})}  \\
    &\left\{\binom{n}{2l-k'-2} \binom{2l-k'-2}{l} \binom{l}{k'+2} L_{k'+1}^2 \left( p_{k'+1}^{-1} - 1 \right) \prod_{i = k'+1}^{l-1}  p_i^{2e_i(\cplx)}  \right\}^{-\frac{1}{2}} \\
    &\left\{\binom{n}{2m-k'-2} \binom{2m-k'-2}{m} \binom{m}{k'+2} M_{k'+1}^2 \left( p_{k'+1}^{-1} - 1 \right) \prod_{i = k'+1}^{m-1} p_i^{2e_i(\mathcal{M})}  \right\}^{-\frac{1}{2}} \\
    &=\lim_{n \to \infty}L_{k'+1}M_{k'+1} \left( p_{k'+1}^{-1} - 1 \right) \prod_{i = k'+1}^{l-1} p_i^{e_i(\cplx) + e_i(\mathcal{M})}\\
   & \left\{  L_{k'+1}^2 \left( p_{k'+1}^{-1} - 1 \right) \prod_{i = k'+1}^{l-1}  p_i^{2e_i(\cplx)}  \prod_{i = k'+1}^{m-1} M_{k'+1}^2 \left( p_{k'+1}^{-1} - 1 \right) \prod_{i = k'+1}^{m-1} p_i^{2e_i(\mathcal{M})}  \right\}^{-\frac{1}{2}}  \\
   &=1.
\end{align*}
\end{proof}

\subsection{Approximation theorem}

Let $\set{\cplx_i}_{i \in [d]}$ be a collection of finite connected simplicial complexes. Let $[m_i]$ be the vertex set of $\cplx_i$ and let $k_i \geq k' + 1$ be its dimension. Let $f_i$ be the function associated to $\cplx_i$ as set out in Equation \ref{equation:f_i}. For any $i \in [d]$, let $\I_i = \binom{[n]}{m_i} \times \set{i}$. For $s = (\phi, i) \in \I_i$ define $Y_{s} = f_i(\mathcal{X}^{(1)}_\phi, \mathcal{X}^{(2)}_\phi, \ldots, \mathcal{X}^{(k_i+1)}_\phi)$ and \[X_s = \sigma_i^{-1}\left\{Y_s - \mu_s \right\},\]
where $\mu_{s} = \EEbracket{Y_{s}}$ and $\sigma_i = \sqrt{\Var(T_{\cplx_i})}$.

Let $W_i = \sum_{s \in \I_i} X_s$ and define a random vector $W = (W_1, W_2, \ldots, W_d) \in \RR^d$. We are interested in the distribution of $W$. To apply the approximation Theorem \ref{theorem:ustat_approx}, we need to show that Assumption \ref{assumption:ho_ustat} holds and find the quantities $\alpha_i$ and $\beta$ from the assumption. The following \tat{proposition} lets us do that.

\begin{proposition}\label{lemma:intersection_of_complexes}
Let $\cplx$ be a connected simplicial complex of dimension $d \geq m$ on the vertex set $[k]$. Consider $s, u \in \binom{[n]}{k}$ such that $|s \cap u| \geq m + 1$. Then there exists $\cplx', \cplx'' \in [\cplx]$ such that $\cplx'[s]$ and $\cplx''[u]$ share at least one $m$-dimensional simplex.
\end{proposition}

Now we can give a normal approximation for subcomplex counts, based on Theorem \ref{theorem:ustat_approx}.
Note that using the definition of $k'$ in this section means that the variables \[\set{\xi^{(1)}_{\alpha}}_{\alpha \in \binom{[n]}{1}}, \set{\xi^{(2)}_{\alpha}}_{\alpha \in \binom{[n]}{2}}, \ldots, \set{\xi^{(k'+1)}_{\alpha}}_{\alpha \in \binom{[n]}{k'+1}}\] are deterministic and equal to 1, which means that when we apply Theorem \ref{theorem:ustat_approx}, we should use $k'+1$ instead of $k'$.

\begin{corollary}\label{theorem:subcomplex_counts}
Let $Z \sim \mvn{0}{\text{\rm Id}_{d \times d}}$ and $\Sigma$ be the covariance matrix of $W$.
\begin{enumerate}
    \item Let $h \in \gr{\mathcal{H}}$. Then \[
        \abs{\EE h(W) - \EE h(\Sigma^{\frac{1}{2}}Z)} \leq \abs{h}_3  B_{\ref{theorem:subcomplex_counts}}  n^{-\frac{k'+2}{2}}.
        \]
        
    \item \gr{Let $m = \max_{i \in [d]} m_i$. Then} \[ \sup _{A \in \mathcal{K}}|\PP(W \in A)-\PP(\Sigma^{\frac{1}{2}}Z \in A)| \leq 2^{\frac{7}{2}} 3^{-\frac{3}{4}}d^{\frac{3}{16}}B_{\ref{theorem:subcomplex_counts}}^{\frac{1}{4}} n^{-\frac{k'+2}{8}}, \]
\end{enumerate}
where
\[
B_{\ref{theorem:subcomplex_counts}} = \frac{4}{3}\frac{(k'+2)^{\frac{3}{2}(k'+2)}}{((k'+2)!)^2} (\sqrt{2}m )^{5m + 2 - 3(k'+2)} d^3 \gamma
\]
and \[\gamma = \left(\prod_{k=k'+1}^{\min_{i \in [d]}(k_i)} p_k^{\min_{i \in [d]}(e_k(\cplx_i))} \right)^{-\frac{1}{2}} \left( 1 - \prod_{k=k'+1}^{\max_{i \in [d]}(k_i)} p_k^{\max_{i \in [d]}(e_k(\cplx_i))} \right) \left( p_{k'+1}^{-1} - 1 \right)^{-\frac{3}{2}}.\]
\end{corollary}

\begin{proof} 
Take $s = (\phi, i), u = (\psi, i) \in \I_i$ such that $|\phi \cap \psi| \geq k'+2$. Using Proposition \ref{lemma:intersection_of_complexes} we get that
\begin{align*}
    \Cov(Y_{s}, Y_{u}) &\geq p_{k'+1}^{2e_{k'+1}(\cplx_i) - 1}\prod_{j=k'+2}^{k_i}p_j^{2e_j(\cplx_i)} + (|[\cplx_i]|^2 - 1)\prod_{j=k'+1}^{k_i}p_j^{2e_j(\cplx_i)} - |[\cplx_i]|^2\prod_{j=k'+1}^{k_i}p_j^{2e_j(\cplx_i)} \\
    &= \prod_{j=k'+1}^{k_i}p_j^{2e_j(\cplx_i)} \left \{ p_{k'+1}^{-1} - 1 \right \} > 0.
\end{align*}

Hence, we see that $\alpha_i = \prod_{j=k'+1}^{k_i}p_j^{2e_j(\cplx_i)} \left \{ p_{k'+1}^{-1} - 1 \right \}$, using the notation of Assumption \ref{assumption:ho_ustat}.

 Recall that $e_k(\cplx_i)$ is the number of $k$-simplices in $\cplx$. For any $i,j,l \in [d]$ and any $s \in \I_i, t \in \I_j, u \in \I_l$ we have using Lemma \ref{lemma:unified_moments} with 
 $\mu_i = \prod_{k=k'+1}^{k_i} p_k^{e_k(\cplx_i)}$,
\begin{align*}
    &\EE \abs{\left\{ Y_{s} - \mu_{s} \right\}\left\{ Y_{t} - \mu_{t} \right\}\left\{ Y_{u} - \mu_{u} \right\}} \leq \\
    &\left\{\prod_{k=k'+1}^{k_i} p_k^{e_k(\cplx_i)} \prod_{k=k'+1}^{k_j} p_k^{e_k(\cplx_j)} \left(1 - \prod_{k=k'+1}^{k_i} p_k^{e_k(\cplx_i)} \right) \left(1 - \prod_{k=k'+1}^{k_j} p_k^{e_k(\cplx_j)} \right) \right\}^{\frac{1}{2}} \\
    & \leq \prod_{k=k'+1}^{\min_{i \in [d]}(k_i)} p_k^{\min_{i \in [d]}(e_k(\cplx_i))} \left( 1 - \prod_{k=k'+1}^{\max_{i \in [d]}(k_i)} p_k^{\max_{i \in [d]}(e_k(\cplx_i))} \right).
\end{align*}

So $\beta = \prod_{k=k'+1}^{\min_{i \in [d]}(k_i)} p_k^{\min_{i \in [d]}(e_k(\cplx_i))} \left( 1 - \prod_{k=k'+1}^{\max_{i \in [d]}(k_i)} p_k^{\max_{i \in [d]}(e_k(\cplx_i))} \right)$.

Now we can apply Theorem \ref{theorem:ustat_approx} and bound the quantity $B_{\ref{theorem:ustat_approx}}$:
\begin{align*}
    B_{\ref{theorem:ustat_approx}} &\leq \frac{4}{3}\frac{(k'+2)^{\frac{3}{2}(k'+2)}}{((k'+2)!)^2} (\sqrt{2}m )^{5m + 2 - 3(k'+2)} d^3 \beta \left \{( p_{k'+1}^{-1} - 1) \prod_{k=k'+1}^{\min_{i \in [d]}(k_i)} p_k^{\min_{i \in [d]}(e_k(\cplx_i))}  \right \}^{-\frac{3}{2}} \\
    & \leq \frac{4}{3}\frac{(k'+2)^{\frac{3}{2}(k'+2)}}{((k'+2)!)^2} (\sqrt{2}m )^{5m + 2 - 3(k'+2)} d^3 \gamma.
\end{align*}
\end{proof}

\begin{remark}
Note that from Corollary \ref{theorem:subcomplex_counts_limiting_cov} it follows that the limiting covariance matrix of $W$, assuming the parameters $p_i$ and the number of vertices $m_i$ are constant for all $i$, is equal to the $d \times d$ matrix with every entry equal to 1. Hence the limiting distribution \tat{in that setting} is a degenerate multivariate normal with covariance matrix having rank 1.

\tat{Also note that if we pick $d$, $p_i \in (0,1)$ for $i > k'$, $k_i$ for all $i$ as well as $e_k(\cplx_i)$ for all $k$ and all $i$ to to be constant, then the bound still goes to zero asymptotically as long as $m = \max_{i \in [d]} m_i$ is of the order $o(\ln^{1 - \epsilon}(n))$ for a fixed $\epsilon \in (0,1)$. By inspecting the inequalities it is possible to vary the parameter setting further since the theorem allows the quantities in the inequalities to depend on each other.}
\end{remark}

\section{Critical simplex counts in lexicographical matchings}\label{section:crit_lexi}

The following lemma is an immediate consequence of Definition \ref{def:lexi_matching} and taken from \cite{temvcinas2021multivariate}.

\begin{lemma}
		\label{proposition:critical_lexi}
		Let $\cplx$ be a simplicial complex endowed with the lexicographical acyclic partial matching, and consider a simplex $t \in \cplx$ with minimal vertex $i \in [n]$. Then, $t$ is matched with
  \begin{enumerate}
      \item one of its co-faces if and only if there exists some $j < i$ for which $t \cup \set{j} \in \cplx$; and,
      \item one of its faces if and only for all $j < i$ we have $(t \setminus \set{i}) \cup \set{j} \notin \cplx$.
  \end{enumerate}  
	\end{lemma}

Let $\mathbf{p} = (p_1, p_2, \ldots, p_{n-1}) \in [0,1]^{n-1}$ and $k' \in \NN$ be such that for any $i \in [k']$ we have $p_i = 1$. Consider the multi-parameter random simplicial complex $\multixnp$ as defined in \cite[Section 2.2]{bobrowski2022random}. For any $k \in [n-1]$ and $\alpha \in \binom{[n]}{k+1}$ let $\xi^{(k+1)}_{\alpha} \sim \ber{p_k}$. 

Fix $s \in \binom{[n]}{k}$. Define the variables $X_s^{+} = \ind{s \text{ matches with its coface given it is a simplex}}$ as well as $X_s^{-} = \ind{s \text{ matches with its face given it is a simplex}}$. The events that the two variables indicate are disjoint. By Lemma \ref{proposition:critical_lexi} we can see that $X_s^{+} = 1 - \prod_{i=1}^{\min(s) - 1}\left(  1 - \prod_{j=k'+1}^{k} \prod_{\substack{\alpha \subseteq s \\ |\alpha| = j}} \xi_{\alpha \cup \set{i}}^{(j+1)} \right)$ and the following: $X_s^{-} = \prod_{i=1}^{\min(s) - 1}\left(  1 - \prod_{j=k'+1}^{k-1} \prod_{\substack{\alpha \subseteq s_{-} \\ |\alpha| = j}} \xi_{\alpha \cup \set{i}}^{(j+1)} \right)$, where $s_{-} \coloneqq s \setminus \set{\min(s)}$. Hence, 
\begin{align*}
    &\ind{s \text{ is a critical simplex}} = \ind{s \in \multixnp}(1 - (X_s^{+} + X_s^{-})) \\
    &= \prod_{j=k'+1}^{k-1} \prod_{\substack{\alpha \subseteq s  \\ |\alpha| = j + 1}} \xi_{\alpha}^{(j+1)} \left[ \prod_{i=1}^{\min(s) - 1}\left(  1 - \prod_{j=k'+1}^{k} \prod_{\substack{\beta \subseteq s \\ |\beta| = j}} \xi_{\beta \cup \set{i}}^{(j+1)} \right) - \prod_{i=1}^{\min(s) - 1}\left(  1 - \prod_{j=k'+1}^{k-1} \prod_{\substack{\beta \subseteq s_{-} \\ |\beta| = j}} \xi_{\beta \cup \set{i}}^{(j+1)} \right)  \right].
\end{align*}

Thus, the random variable of interest is \tat{the number of $k$-simplices that are critical under the lexicographical matching,}
\begin{equation}\label{equation:def_lexi_crit}
T_k = \sum_{s \in \binom{[n]}{k}}\prod_{j=k'+1}^{k-1} \prod_{\substack{\alpha \subseteq s  \\ |\alpha| = j + 1}} \xi_{\alpha}^{(j+1)} \left[ \prod_{i=1}^{\min(s) - 1}\left(  1 - \prod_{j=k'+1}^{k} \prod_{\substack{\beta \subseteq s  \\ |\beta| = j}} \xi_{\beta \cup \set{i}}^{(j+1)} \right) - \prod_{i=1}^{\min(s) - 1}\left(  1 - \prod_{j=k'+1}^{k-1} \prod_{\substack{\beta \subseteq s_{-} \\ |\beta| = j}} \xi_{\beta \cup \set{i}}^{(j+1)} \right)  \right].
\end{equation}
Note that this random variable does not fit into the framework of generalised $U$-statistics \tat{as in Definition \ref{definition:u_statistic_order_k}} because the summands depend not only on the variables that are indexed by the subset $s$.

\subsection{Moments and the limiting covariance} \label{subsec:criticalmoments}

\begin{lemma}\label{lemma:critlexi_mean}
For $k' + 1 \leq k \leq n - 1$ we have:
\[\binom{n-2}{k} \prod_{i=k'+1}^k p_i^{\binom{k+1}{i+1} + \binom{k}{i}} \left( 1 - \prod_{i=k'+1}^{k+1}p_i^{\binom{k}{i-1}} \right) \leq \EEbracket{T_{k+1}} \leq \binom{n-1}{k} \prod_{i=k'+1}^{k+1}p_i^{\binom{k+1}{i+1} - \binom{k+1}{i}} \left( 1 - \prod_{i=k'+1}^{k} p_i^{\binom{k}{i-1}} \right).\]
\end{lemma}

\begin{proof}
\begin{align*}
    \EEbracket{T_{k+1}} &= \sum_{l=1}^{n-k} \sum_{\substack{s \in \binom{[n]}{k+1} \\ \min(s) = l}} \prod_{i=k'+1}^k p_i^{\binom{k+1}{i+1}} \left[ \left( 1 - \prod_{i=k'+1}^{k+1}p_i^{\binom{k+1}{i}} \right)^{l-1} - \left( 1 - \prod_{i=k'+1}^{k}p_i^{\binom{k}{i}} \right)^{l-1}  \right] \\
    &= \sum_{l=1}^{n-k} \binom{n-l}{k} \prod_{i=k'+1}^k p_i^{\binom{k+1}{i+1}} \left[ \left( 1 - \prod_{i=k'+1}^{k+1}p_i^{\binom{k+1}{i}} \right)^{l-1} - \left( 1 - \prod_{i=k'+1}^{k}p_i^{\binom{k}{i}} \right)^{l-1}  \right] \\
    &\leq \binom{n-1}{k} \prod_{i=k'+1}^k p_i^{\binom{k+1}{i+1}} \sum_{l=0}^{\infty}  \left[ \left( 1 - \prod_{i=k'+1}^{k+1}p_i^{\binom{k+1}{i}} \right)^{l} - \left( 1 - \prod_{i=k'+1}^{k}p_i^{\binom{k}{i}} \right)^{l}  \right] \\
    &= \binom{n-1}{k} \prod_{i=k'+1}^{k+1}p_i^{\binom{k+1}{i+1} - \binom{k+1}{i}} \left( 1 - \prod_{i=k'+1}^{k} p_i^{\binom{k}{i-1}} \right).
\end{align*}

For the lower bound, since all the summands are non-negative, we bound \tat{the expectation} by the $l=2$ term:
\begin{align*}
    \EEbracket{T_{k+1}} &\geq \binom{n-2}{k} \prod_{i=k'+1}^k p_i^{\binom{k+1}{i+1}} \left[ \left( 1 - \prod_{i=k'+1}^{k+1}p_i^{\binom{k+1}{i}} \right) - \left( 1 - \prod_{i=k'+1}^{k}p_i^{\binom{k}{i}} \right)  \right] \\
    & = \binom{n-2}{k} \prod_{i=k'+1}^k p_i^{\binom{k+1}{i+1} + \binom{k}{i}} \left( 1 - \prod_{i=k'+1}^{k+1}p_i^{\binom{k}{i-1}} \right).
\end{align*}
\end{proof}

\tat{The following two lemmas follow similar arguments as in Lemma \ref{lemma:critlexi_mean} but are more involved. The proofs are long calculations which are not particularly insightful and so they are deferred to the Appendix \ref{section:appendix}.}

\begin{lemma}\label{lemma:lower_bound_var_crit_lexi}
For a fixed integer $k' + 1 \leq k \leq n-1$ and a fixed sequence of probabilities $p_i \in [0,1]$ there is a constant $C > 0$ independent of $n$ and a natural number $N_{p,k}$ such that for any $n \geq N_{p,k}$:
\[\Var(T_{k+1}) \geq C n^{2k-k'}.\]
\end{lemma}

\begin{lemma}\label{theorem:crit_lexi_limiting_cov}
Let $\sigma_{i+k'}^{2} = \Var(T_{i+k'+1})$ for any $i \geq 1$. Let $\Sigma$ be a $d \times d$ symmetric matrix such that $\Sigma_{i,j} = \sigma_{i+k'}^{-1} \sigma_{j + k'}^{-1} \Cov(T_{i+k'+1}, T_{j+k'+1})$. Then assuming that for all $l$ the parameters $p_l$ do not depend on $n$, for any fixed $i < j$ we have:
\[\lim_{n \to \infty} \Sigma_{i,j} =  \frac{\rho(k+1)\rho(r+1) \sqrt{(2 - \rho(k+1))(2 - \rho(r+1))(2 - \rho(k+1)\rho(k'+1)^{-1})(2 - \rho(r+1)\rho(k'+1)^{-1})}}{(\rho(k+1) +\rho(r+1) - \rho(k+1)\rho(r+1)\rho(k'+1)^{-1})(\rho(k+1) +\rho(r+1) - \rho(k+1)\rho(r+1))},\]
where $k'+i = k$, $k'+j = r$, and $\rho(a) \coloneqq \prod_{i=k'+1}^{a}p_i^{\binom{a}{i}}$.
\end{lemma}

\gr{\begin{remark} The square root in Lemma \ref{theorem:crit_lexi_limiting_cov} is understood to be the positive square root. As $0 \le \rho(a) \le 1$ the limiting expression in Lemma \ref{theorem:crit_lexi_limiting_cov} is always non-negative. Moreover we can bound the limiting covariance from above by $\frac{4\rho(k+1)\rho(r+1)}{(\rho(k+1) + \rho(r+1))^2}$. Note that $\frac{4\rho(k+1)\rho(r+1)}{(\rho(k+1) + \rho(r+1))^2} < 1$ as long as $\rho(k+1) \neq \rho(r+1)$.  In particular, if $p_k \in (0,1)$ for all $k > k'$, then  $\rho(k+1) \neq \rho(r+1)$ and so the limiting covariance is strictly smaller than 1.
\end{remark}}

\subsection{Approximation theorem}

For $i \in [d]$, recall a random variable counting $i$-simplices in $\multixnp$ that are critical under the lexicographical matching, as given in \eqref{equation:def_lexi_crit}. We write for the $i$-th index set $\I_i \coloneqq \binom{[n]}{i+1} \times \set{i}$. For $s = (\phi,i) \in \I_i$ we \tat{have} 
\[\mu_{s} = \prod_{j=k'+1}^k p_j^{\binom{k+1}{j+1}} \left[ \left( 1 - \prod_{i=k'+1}^{k+1}p_i^{\binom{k+1}{i}} \right)^{\min(\phi)-1} - \left( 1 - \prod_{i=k'+1}^{k}p_i^{\binom{k}{i}} \right)^{\min(\phi)-1}  \right]\] 
and $\sigma_i = \sqrt{\Var(T_{i+1})}$. \tat{Note that $\EEbracket{T_{i+1}} = \sum_{s \in \I_i} \mu_s$.} Let 
\[X_{s} = \sigma_i^{-1}\left\{ \prod_{j=k'+1}^{k-1} \prod_{\substack{\alpha \subseteq \phi  \\ |\alpha| = j + 1}} \xi_{\alpha}^{(j+1)} \left[ \prod_{i=1}^{\min(\phi) - 1}\left(  1 - \prod_{j=k'+1}^{k} \prod_{\substack{\alpha \subseteq \phi \\ |\alpha| = j}} \xi_{\alpha \cup \set{i}}^{(j+1)} \right) - \prod_{i=1}^{\min(\phi) - 1}\left(  1 - \prod_{j=k'+1}^{k-1} \prod_{\substack{\alpha \subseteq \phi_{-} \\ |\alpha| = j}} \xi_{\alpha \cup \set{i}}^{(j+1)} \right)  \right]  - \mu_{s}\right\}.\]

Let $W_i = \sum_{s \in \I_i} X_{s}$ and $W = (W_{k'+1}, W_{k'+2}, \ldots, W_{k'+d}) \in \RR^d$. For bounds that asymptotically go to zero for this example, we use Theorem \ref{theorem:mvn_dissociated_decomp_approx} directly:  the uniform bounds from Corollary \ref{corollary:mvn_dissociated_decomp} are not fine enough.

\begin{theorem}\label{theorem:crit_lexi_approx}
Let $Z \sim \mvn{0}{\text{\rm Id}_{d \times d}}$ and $\Sigma$ be the covariance matrix of $W$.
\begin{enumerate}
    \item Let $h \in \gr{\mathcal{H}}$. Then there is a constant $B_{\ref{theorem:crit_lexi_approx}.1} >0$ independent of $n$ and a natural number $N_{\ref{theorem:crit_lexi_approx}.1}$ such that for any $n \geq N_{\ref{theorem:crit_lexi_approx}.1}$ we have\[\abs{\EE h(W) - \EE h(\Sigma^{\frac{1}{2}}Z)} \leq B_{\ref{theorem:crit_lexi_approx}.1} \abs{h}_3 n^{-1-\frac{k'}{2}}.\]
    \item There is a constant $B_{\ref{theorem:crit_lexi_approx}.2} >0$ independent of $n$ and a natural number $N_{\ref{theorem:crit_lexi_approx}.2}$ such that for any $n \geq N_{\ref{theorem:crit_lexi_approx}.2}$ we have\[\sup _{A \in \mathcal{K}}|\PP(W \in A) -\PP(\Sigma^{\frac{1}{2}}Z \in A)| \leq B_{\ref{theorem:crit_lexi_approx}.2}n^{-\frac{1}{4} - \frac{k'}{8}}.\]
\end{enumerate}
\end{theorem}

\begin{proof}  It is clear that $W$ satisfies the conditions of Theorem \ref{theorem:mvn_dissociated_decomp_approx} 
for any $s = (\phi,i) \in \I_i$ setting
\[\D_j(s) = \setc{(\psi, j) \in \I_j}{ |\phi \cap \psi| \geq k' + 1 }.\]
We apply Theorem \ref{theorem:mvn_dissociated_decomp_approx}. We write $C$ for an unspecified positive constant that does not depend on $n$. For the bounds on the quantity $B_{\ref{theorem:mvn_dissociated_decomp_approx}}$ we use Lemma \ref{lemma:unified_moments} and Lemma \ref{lemma:lower_bound_var_crit_lexi} as well as the fact that for fixed $b \in [n], i,j \in [k'+1, k'+d], s \in \I_i$ we have $|\D_j(s)| \leq Cn^{j-k'}$ and $|\setc{(\psi, j) \in \D_j(s)}{\min(\psi) = b}| \leq Cn^{j-k'-1}$.  Also, we assume here that $n$ is large enough for the bound in Lemma \ref{lemma:lower_bound_var_crit_lexi} to apply. \tat{The existence of $B_{\ref{theorem:crit_lexi_approx}.1}$ and $B_{\ref{theorem:crit_lexi_approx}.2}$ is proved by applying Theorem \ref{theorem:mvn_dissociated_decomp_approx} and bounding the quantity $B_{\ref{theorem:mvn_dissociated_decomp_approx}}$ from the theorem.}
\begin{align*}
&B_{\ref{theorem:mvn_dissociated_decomp_approx}} \leq \frac{1}{3} \sum_{i,j,k=1}^d \sum_{a=1}^{n-i} \sum_{\substack{\phi \in \binom{[n]}{i+1} \\ \min(\phi) = a}} \sum_{b=1}^{n-j} \sum_{\substack{(\psi,j) \in \D_j((\phi,i)) \\ \min(\psi) = b}} \Bigg\{ \sum_{r \in \D_k((\phi,i))} \frac{3}{2}  (\sigma_i \sigma_j \sigma_k)^{-1} \left\{ \mu(a)\mu(b)(1-\mu(a))(1-\mu(b)) \right\}^{\frac{1}{2}}  \\
&+ \sum_{r \in \D_k((\psi,j))} (\sigma_i \sigma_j \sigma_k)^{-1} \left\{ \mu(a)\mu(b)(1-\mu(a))(1-\mu(b)) \right\}^{\frac{1}{2}} \Bigg\}\\
&\leq  \sum_{i,j,k=1}^d \sum_{a=1}^{n-i} \sum_{b=1}^{n-j} Cn^{i+j+k-2k'-1}n^{\frac{3}{2}k'-i-j-k} \left\{ (1-\rho(i+1))^{a - 1}(1-\rho(j+1))^{b - 1} + (1-\rho(i))^{a - 1}(1-\rho(j))^{b - 1} \right\}^{\frac{1}{2}} 
\\
&\leq Cn^{-1-\frac{k'}{2}} d^3 \sum_{a=1}^{\infty} \sum_{b=1}^{\infty} \left[\left\{ (1-\rho(i+1))^{a - 1}(1-\rho(j+1))^{b - 1} \right\}^{\frac{1}{2}} + \left\{ (1-\rho(i))^{a - 1}(1-\rho(j))^{b - 1} \right\}^{\frac{1}{2}} \right] \leq  C  n^{-1-\frac{k'}{2}}. \\
\end{align*}
\end{proof}
\subsection*{Acknowledgements}
	TT acknowledges funding from EPSRC studentship 2275810. VN and GR are supported by the EPSRC grant EP/R018472/1. GR is also funded in part by the EPSRC grant EP/T018445/1. The authors would like to thank Christina Goldschmidt, Heather Harrington,  \gr{Adrian R{\"o}llin and Fang Xiao} for helpful discussions.

\printbibliography

\newpage
\begin{appendix}
		\section{Proof of Lemmas \ref{lemma:lower_bound_var_crit_lexi} and \ref{theorem:crit_lexi_limiting_cov}}\label{section:appendix}

  \begin{lemma}\label{lemma:var_crit_lexi}
For any integer $k'+1 \leq k \leq n - 1$ we have: \[\Var\{ T_{k+1} \} = V_1 + V_2 + V_3 + V_4,\]

where

\begin{align*}
&V_1 = 2\sum_{l < m}^{n-k} \sum_{j=k'+1}^k \sum_{q=1}^{\min(k+1, m-l)} \binom{n-m}{2k+1-j-q} \binom{2k+1-j-q}{k} \binom{k}{j-1} \binom{m-l+1}{q-1} \\
&\Big\{ \rho^+(k+1)^2\rho^+(j)^{-1} \tau(l,m,q,k,j-1) \\
&\left[ \left(1-2\rho(k) + \rho(k)^2\rho(j-1)^{-1}\right)^{l-1} - \left(1-\rho(k+1)-\rho(k) + \rho(k+1)\rho(k)\rho(j-1)^{-1}\right)^{l-1}\right] \\
&+ \rho^+(k+1)^2\rho^+(j)^{-1} \tau(l,m,q,k+1,j) \\
&\left[ \left(1-2\rho(k+1) + \rho(k+1)^2\rho(j)^{-1} \right)^{l-1} - \left(1-\rho(k+1)-\rho(k) + \rho(k)\rho(k+1)\rho(j)^{-1}\right)^{l-1}\right] - \mu(l)\mu(m) \Big\} ;\\
&V_2 = 2\sum_{l < m}^{n-k} \sum_{j=k'+1}^k \sum_{q=1}^{\min(k+1, m-l)} \binom{n-m}{2k+1-j-q} \binom{2k+1-j-q}{k} \binom{k}{j} \binom{m-l+1}{q-1} \\
&\Big\{ \rho^+(k+1)^2\rho^+(j)^{-1} \tau(l,m,q,k,j) \\
&\left[ \left(1-2\rho(k) + \rho(k)^2\rho(j)^{-1}\right)^{l-1} - \left(1-\rho(k+1)-\rho(k) + \rho(k+1)\rho(k)\rho(j)^{-1}\right)^{l-1}\right] \\
&+ \rho^+(k+1)^2\rho^+(j)^{-1} \tau(l,m,q,k+1,j) \\
&\left[ \left(1-2\rho(k+1) + \rho(k+1)^2\rho(j)^{-1} \right)^{l-1} - \left(1-\rho(k+1)-\rho(k) + \rho(k)\rho(k+1)\rho(j)^{-1}\right)^{l-1}\right] - \mu(l)\mu(m) \Big\} ;\\
&V_3 = \sum_{l=1}^{n-k} \sum_{j=k'+1}^k  \binom{n-l}{2k+1-j} \binom{2k+1-j}{k} \binom{k}{j-1} \Big\{ \rho^{+}(k+1)^2\rho^{+}(j)^{-1} \\
&\Big[ \left(1-2\rho(k+1) + \rho(k+1)^2\rho(j)^{-1} \right)^{l-1} + \left(1-2\rho(k) + \rho(k)^2\rho(j-1)^{-1} \right)^{l-1} \\
&- 2\left(1 -\rho(k) -\rho(k+1) + \rho(k)\rho(k+1)\rho(j-1)^{-1}\right)^{l-1} \Big] -  \mu(l)^2\Big\}; \\
&V_4 = \sum_{l=1}^{n-k} \binom{n-l}{k} \left\{\mu(l) - \mu(l)^2\right\}.
\end{align*}

Here we have used the following abbreviations:
\begin{align*}
&\rho(a) \coloneqq \prod_{i=k'+1}^{a}p_i^{\binom{a}{i}}; &&\tau(l,m,q,a,b) \coloneqq \left(1-\rho(a)\right)^{m-l-q} \left(1-\rho(a)\rho(b)^{-1}\right)^{q};\\
&\rho^{+}(a) \coloneqq \prod_{i=k'+1}^{a}p_i^{\binom{a}{i+1}}; &&\mu(a) \coloneqq \rho^{+}(k+1)\left(\left(1-\rho(k+1)\right)^{a-1} - \left(1-\rho(k)\right)^{a-1}\right).
\end{align*}

Also, the notation $ \sum_{l < m}^{n-k}$ stands for $ \sum_{l=1}^{n-k-1} \sum_{m= l+1}^{n-k}$.
\end{lemma}

\begin{proof}[Proof of Lemma \ref{lemma:var_crit_lexi}]
Here, we adapt and generalise the proof of \cite[Lemma 34]{temvcinas2021multivariate}. For $s \in \binom{[n]}{k+1}$ recall that $s_{-} = s \setminus \set{\min(s)}$. We write:
\begin{align*}
    Y^+_s &= \prod_{i=1}^{\min(s) - 1}\left(  1 - \prod_{j=k'+1}^{k+1} \prod_{\substack{\alpha \subseteq s \\ |\alpha| = j}} \xi_{\alpha \cup \set{i}}^{(j+1)} \right),
    &Y_s^{-} &= \prod_{i=1}^{\min(s) - 1}\left(  1 - \prod_{j=k'+1}^{k} \prod_{\substack{\alpha \subseteq s_{-}  \\ |\alpha| = j }} \xi_{\alpha \cup \set{i}}^{(j+1)} \right),\\
    Z_s &=  \prod_{j=k'+1}^{k} \prod_{\substack{\alpha \subseteq s  \\ |\alpha| = j + 1}} \xi_{\alpha}^{(j+1)},
    &Y_s &= Y_s^{+} - Y_s^{-}.
\end{align*}
Then $Z_s$ and $Y_s$ are independent and $T_{k+1} = \sum_{s \in \binom{[n]}{k+1}} Z_s Y_s$. Consider the variance: 
\begin{align}\label{varexapansion}
    &\Var(T_{k+1}) = \sum_{s \in \binom{[n]}{k+1}} \Var(Z_sY_s) + \sum_{\substack{s \neq t \in \binom{[n]}{k+1} \\ \min(s) \neq \min(t)}} \Cov(Z_sY_s, Z_tY_t) + \sum_{\substack{s \neq t \in \binom{[n]}{k+1} \\ \min(s) = \min(t)}} \Cov(Z_sY_s, Z_tY_t). 
\end{align}
For the first term in \eqref{varexapansion}, writing $\min(s) = l$, we get: \[ \PP(Z_sY_s=1)= \mu(l) = \prod_{i=k'+1}^k p_i^{\binom{k+1}{i+1}} \left[ \left( 1 - \prod_{i=k'+1}^{k+1}p_i^{\binom{k+1}{i}} \right)^{l-1} - \left( 1 - \prod_{i=k'+1}^{k}p_i^{\binom{k}{i}} \right)^{l-1}  \right].\] Now we see:
\begin{align*}
&\sum_{s \in \binom{[n]}{k+1}} \Var(Z_sY_s) = \sum_{l=1}^n \sum_{\substack{s \in \binom{[n]}{k+1} \\ \min(s) = l}} \left( \EEbracket{(Z_sY_s)^2} - \EEbracket{Z_sY_s}^2 \right) \\
&= \sum_{l=1}^{n-k}  \binom{n-l}{k} \left( \PP(Z_sY_s=1) - \PP(Z_sY_s=1)^2 \right) = \sum_{l=1}^{n-k} \binom{n-l}{k} \left\{ \mu(l) - \mu(l)^2 \right\} = V_4.
\end{align*}

Now consider the covariance terms in \eqref{varexapansion}, the expansion of the variance. Note that for any $s, t \in \binom{[n]}{k+1}$ if $|s \cap t| < k' + 1$, then the variables $Z_s Y_s$ and $Z_t Y_t$ can be written as functions of two disjoint sets of independent random variables for the form $\xi_{\alpha}^{(l)}$ for $\alpha \in \binom{[n]}{l}$ and hence have zero covariance.

Fix $s, t \in \binom{[n]}{k+1}$ and assume $k' + 1 \leq |s \cap t| \leq k$. Note that because $|s \cap t| \neq k + 1$, we have $s \neq t$.
There are $2\binom{k+1}{i} - \binom{|s \cap t|}{i}$ distinct simplices of size $i$ in $s$ and $t$ combined and hence $\PP(Z_s Z_t = 1) = \prod_{i=k'+1}^k p_i^{2\binom{k+1}{i+1} - \binom{|s \cap t|}{i}}$. Also, $Y_sY_t = Y_s^+Y_t^+ + Y_s^-Y_t^- - Y_s^+Y_t^- - Y_s^-Y_t^+$. For the rest of the proof when calculating probabilities we  assume w.l.o.g.\,that $\min(s) \leq \min(t)$. Then we have for $Y_s^+ Y_t^+$:

\begin{align*}
Y_s^+Y_t^+ &= \prod_{i=1}^{\min(t) - 1}\left(  1 - \prod_{j=k'+1}^{k+1} \prod_{\substack{\alpha \subseteq t \\ |\alpha| = j}} \xi_{\alpha \cup \set{i}}^{(j+1)} \right) \prod_{i=1}^{\min(s) - 1}\left(  1 - \prod_{j=k'+1}^{k+1} \prod_{\substack{\alpha \subseteq s \\ |\alpha| = j}} \xi_{\alpha \cup \set{i}}^{(j+1)} \right) \\
&=\prod_{i=1}^{\min(s) - 1}\left\{ \left(  1 - \prod_{j=k'+1}^{k+1} \prod_{\substack{\alpha \subseteq s \\ |\alpha| = j}} \xi_{\alpha \cup \set{i}}^{(j+1)} \right)\left(  1 - \prod_{j=k'+1}^{k+1} \prod_{\substack{\alpha \subseteq t \\ |\alpha| = j}} \xi_{\alpha \cup \set{i}}^{(j+1)} \right) \right\}  \prod_{\substack{i=\min(s)\\i \in s}}^{\min(t) - 1}\left(  1 - \prod_{j=k'+1}^{k+1} \prod_{\substack{\alpha \subseteq t \\ |\alpha| = j}} \xi_{\alpha \cup \set{i}}^{(j+1)} \right) \\
&\prod_{\substack{i=\min(s)\\ i \notin s}}^{\min(t) - 1}\left(  1 - \prod_{j=k'+1}^{k+1} \prod_{\substack{\alpha \subseteq t \\ |\alpha| = j}} \xi_{\alpha \cup \set{i}}^{(j+1)} \right).
\end{align*}

Fix $i \in [\min(s)-1]$. Then with $\neg$ denoting the complement
\begin{align*}
& \PP \left[ \left(  1 - \prod_{j=k'+1}^{k+1} \prod_{\substack{\alpha \subseteq s \\ |\alpha| = j}} \xi_{\alpha \cup \set{i}}^{(j+1)} \right)\left(  1 - \prod_{j=k'+1}^{k+1} \prod_{\substack{\alpha \subseteq t \\ |\alpha| = j}} \xi_{\alpha \cup \set{i}}^{(j+1)} \right) = 1 \right] \\
&= \PP \left[ \neg\left(\prod_{j=k'+1}^{k+1} \prod_{\substack{\alpha \subseteq s \\ |\alpha| = j}} \xi_{\alpha \cup \set{i}}^{(j+1)} = 1 \cup \prod_{j=k'+1}^{k+1} \prod_{\substack{\alpha \subseteq t \\ |\alpha| = j}} \xi_{\alpha \cup \set{i}}^{(j+1)} = 1 \right) \right] =1-2\prod_{i=k'+1}^{k+1}p_i^{\binom{k+1}{i}} + \prod_{i=k'+1}^{k+1}p_i^{2\binom{k+1}{i} - \binom{|s \cap t|}{i}}.
\end{align*}
Here  $\cup$  indicates that either or both of these events may happen.

Moreover, the variable $\prod_{i=1}^{\min(s) - 1}\left(  1 - \prod_{j=k'+1}^{k+1} \prod_{\substack{\alpha \subseteq s \\ |\alpha| = j}} \xi_{\alpha \cup \set{i}}^{(j+1)} \right)\left(  1 - \prod_{j=k'+1}^{k+1} \prod_{\substack{\alpha \subseteq t \\ |\alpha| = j}} \xi_{\alpha \cup \set{i}}^{(j+1)} \right)$ and the variable $\prod_{\substack{i=\min(s)\\ i \notin s}}^{\min(t) - 1}\left(  1 - \prod_{j=k'+1}^{k+1} \prod_{\substack{\alpha \subseteq t \\ |\alpha| = j}} \xi_{\alpha \cup \set{i}}^{(j+1)} \right)$ are independent of  $Z_s Z_t$. Recall the notation $[a,b] = \set{a, a+1, \ldots, b}$ for two positive integers $a \leq b$. Setting $q \coloneqq |s \cap [\min(s), \min(t) - 1]|$, 
\begin{align*}
&\PP\left(Y_s^+ Y_t^+ = 1 | Z_s Z_t = 1\right) \\ 
= &\PP\left(\prod_{i=1}^{\min(s) - 1}\left(  1 - \prod_{j=k'+1}^{k+1} \prod_{\substack{\alpha \subseteq s \\ |\alpha| = j}} \xi_{\alpha \cup \set{i}}^{(j+1)} \right)\left(  1 - \prod_{j=k'+1}^{k+1} \prod_{\substack{\alpha \subseteq t \\ |\alpha| = j}} \xi_{\alpha \cup \set{i}}^{(j+1)} \right) =1 \right) \\ &\PP\left(\prod_{\substack{i=\min(s)\\ i \notin s}}^{\min(t) - 1}\left(  1 - \prod_{j=k'+1}^{k+1} \prod_{\substack{\alpha \subseteq t \\ |\alpha| = j}} \xi_{\alpha \cup \set{i}}^{(j+1)} \right) =1 \right)\PP\left(\prod_{\substack{i=\min(s)\\ i \in s}}^{\min(t) - 1}\left(  1 - \prod_{j=k'+1}^{k+1} \prod_{\substack{\alpha \subseteq t \\ |\alpha| = j}} \xi_{\alpha \cup \set{i}}^{(j+1)} \right) =1 \Bigg| Z_s Z_t = 1 \right)\\
=&\left(1-2\prod_{i=k'+1}^{k+1}p_i^{\binom{k+1}{i}} + \prod_{i=k'+1}^{k+1}p_i^{2\binom{k+1}{i} - \binom{|s \cap t|}{i}}\right)^{\min(s) - 1}\left(1-\prod_{i=k'+1}^{k+1}p_i^{\binom{k+1}{i}}\right)^{\min(t) - \min(s)-q} \left(1-\prod_{i=k'+1}^{k+1}p_i^{\binom{k+1}{i} - \binom{|s \cap t|}{i}}\right)^{q}.
\end{align*}

This strategy of splitting the product $Y_s^+Y_t^+$ into three products of independent variables, only one of which is dependent on $Z_sZ_t$ works exactly in the same way for the variables $Y_s^-Y_t^+$, $Y_s^+Y_t^-$, $Y_s^-Y_t^-$. We write $l = \min(s)$, $m = \min(t)$. Also, we set:
\begin{align*} \pi(l,m,a, b, d_1, d_2, q) \coloneqq  &
\left(1-\prod_{i=k'+1}^{a}p_i^{\binom{a}{i}} - \prod_{i=k'+1}^{b}p_i^{\binom{b}{i}} + \prod_{i=k'+1}^{\max(a,b)}p_i^{\binom{a}{i} + \binom{b}{i} - \binom{d_1}{i}}\right)^{l-1}\left(1-\prod_{i=k'+1}^{a}p_i^{\binom{a}{i}}\right)^{m-l-q}\\
&\left(1-\prod_{i=k'+1}^{a}p_i^{\binom{a}{i} - \binom{d_2}{i}}\right)^q.
\end{align*} 

Using the described strategy we get:
\begin{align*}
    &\PP\left(Y_s^+ Y_t^+ = 1 | Z_s Z_t = 1\right) = \pi(l,m,k+1, k+1, |s \cap t|, |s \cap t|, q)\\
    &\PP\left(Y_s^- Y_t^- = 1 | Z_s Z_t = 1\right) = \pi(l,m,k, k, |s_{-} \cap t_{-}|, |s \cap t_{-}|, q)\\
    &\PP\left(Y_s^+ Y_t^- = 1 | Z_s Z_t = 1\right) = \pi(l,m,k, k+1, |s \cap t_{-}|, |s \cap t_{-}|, q)\\
    &\PP\left(Y_s^- Y_t^+ = 1 | Z_s Z_t = 1\right) = \pi(l,m,k+1,k,|s_{-} \cap t|,|s \cap t|,q).
\end{align*}

Now we are ready to calculate the covariance:
\begin{align*}
&\Cov(Z_sY_s, Z_tY_t) = \EEbracket{Z_sZ_tY_s^+Y_t^+} + \EEbracket{Z_sZ_tY_s^-Y_t^-} - \EEbracket{Z_sZ_tY_s^+Y_t^-} - \EEbracket{Z_sZ_tY_s^-Y_t^+} \\
&- \EEbracket{Z_sY_s}\EEbracket{Z_tY_t}\\
&= \PP(Z_s Z_t = 1)\Big\{\PP\left(Y^+_s Y^+_t = 1 | Z_s Z_t = 1\right) + \PP\left(Y_s^- Y_t^- = 1 | Z_s Z_t = 1\right) \\
&- \PP\left(Y^+_s Y^-_t = 1 | Z_s Z_t = 1\right) - \PP\left(Y^-_s Y^+_t = 1 | Z_s Z_t = 1\right) \Big\} - \PP(Z_s Y_s = 1)\PP(Z_t Y_t = 1)\\ 
&= \prod_{i=k'+1}^k p_i^{2\binom{k+1}{i+1} - \binom{|s \cap t|}{i+1}} (\pi(l,m,k+1, k+1, |s \cap t|, |s \cap t|, q) + \pi(l,m,k, k, |s_{-} \cap t_{-}|, |s \cap t_{-}|, q) \\
&- \pi(l,m,k, k+1, |s \cap t_{-}|, |s \cap t_{-}|, q) - \pi(l,m,k+1,k,|s_{-} \cap t|,|s \cap t|,q)) -\mu(l)\mu(m).
\end{align*}

Next we consider the two covariance sums in \eqref{varexapansion} separately. First assume that $\min(s) \ne \min (t)$. Given $l,m \in [n-k]$, $j \in [k]$, and $q \in [\min(k+1, |m-l|)]$ define the set 
\[\Gamma_{k+1}(l,m,j,q) = \setc{(s,t)}{s,t \in \binom{[n]}{k+1}, \min(s) = l, \min(t) = m, |s \cap t| = j, \max(q_{s,t}, q_{t,s}) = q}\]
as well as
\[\Gamma^+_{k+1}(l,m,j,q) = \setc{(s,t) \in \Gamma_{k+1}(l,m,j,q)}{\min(t) \in s}\]
and 
\[\Gamma^-_{k+1}(l,m,j,q) = \setc{(s,t) \in \Gamma_{k+1}(l,m,j,q)}{\min(t) \notin s}.\]
Here $q_{s,t} = |s \cap [\min(s), \min(t) - 1]|, q_{t,s} = |t \cap [\min(t), \min(s) - 1]|$. From the proof of Lemma 4.3 in \cite{temvcinas2021multivariate} we know that: 
\[
|\Gamma^+_{k+1}(l,m,j,q)| = \binom{n-m}{2k+1-j-q} \binom{2k+1-j-q}{k} \binom{k}{j-1} \binom{m-l+1}{q-1}; 
\]
 and 
\[
|\Gamma^-_{k+1}(l,m,j,q)| = \binom{n-m}{2k+1-j-q} \binom{2k+1-j-q}{k} \binom{k}{j} \binom{m-l+1}{q-1}. 
\]

Now using the covariance expression we have just derived, we get
\begin{align*}
&\sum_{\substack{s \neq t \in \binom{[n]}{k+1} \\ \min(s) \neq \min(t)}} \Cov(Z_sY_s, Z_tY_t)\\
=  &\sum_{l=1}^{n-k} \sum_{m=i+1}^{n-k} \sum_{j=k'+1}^k \sum_{q=1}^{\min(k+1, m-l)} \sum_{(s,t) \in \Gamma^+_{k+1}(l,m,j,q)} \Cov(Z_sY_s, Z_tY_t) \\
&+ \sum_{l=1}^{n-k} \sum_{m=i+1}^{n-k} \sum_{j=k'+1}^k \sum_{q=1}^{\min(k+1, m-l)} \sum_{(s,t) \in \Gamma^-_{k+1}(l,m,j,q)} \Cov(Z_sY_s, Z_tY_t) \\
&+ \sum_{m=1}^{n-k} \sum_{l=j+1}^{n-k} \sum_{j=k'+1}^k \sum_{q=1}^{\min(k+1, l-m)} \sum_{(s,t) \in \Gamma^+_{k+1}(m,l,j,q)} \Cov(Z_sY_s, Z_tY_t) \\
&+ \sum_{m=1}^{n-k} \sum_{l=j+1}^{n-k} \sum_{j=k'+1}^k \sum_{q=1}^{\min(k+1, l-m)} \sum_{(s,t) \in \Gamma^-_{k+1}(m,l,j,q)} \Cov(Z_sY_s, Z_tY_t)\\
= &\sum_{l=1}^{n-k} \sum_{m=i+1}^{n-k} \sum_{j=k'+1}^k \sum_{q=1}^{\min(k+1, m-l)} |\Gamma^+_{k+1}(l,m,j,q)| \Big\{ \prod_{i=k'+1}^k p_i^{2\binom{k+1}{i+1} - \binom{j}{i+1}} (\pi(l,m,k+1, k+1, j, j, q) \\
&+ \pi(l,m,k, k, j-1, j-1, q) - \pi(l,m,k, k+1, j-1, j-1, q) - \pi(l,m,k+1,k,j,j,q)) -\mu(l)\mu(m) \Big\}\\
&+ \sum_{l=1}^{n-k} \sum_{m=i+1}^{n-k} \sum_{j=k'+1}^k \sum_{q=1}^{\min(k+1, m-l)} |\Gamma^-_{k+1}(l,m,j,q)| \Big\{ \prod_{i=k'+1}^k p_i^{2\binom{k+1}{i+1} - \binom{j}{i+1}} (\pi(l,m,k+1, k+1, j, j, q) \\
&+ \pi(l,m,k, k, j,j, q) - \pi(l,m,k, k+1, j, j, q)  - \pi(l,m,k+1,k,j,j,q)) -\mu(l)\mu(m) \Big\}\\
&+ \sum_{m=1}^{n-k} \sum_{l=j+1}^{n-k} \sum_{j=k'+1}^k \sum_{q=1}^{\min(k+1, l-m)} |\Gamma^+_{k+1}(m,l,j,q)| \Big\{ \prod_{i=k'+1}^k p_i^{2\binom{k+1}{i+1} - \binom{j}{i+1}} (\pi(m,l,k+1, k+1, j,j, q) \\
&+ \pi(m,l,k, k, j-1, j-1, q) - \pi(m,l,k, k+1, j-1, j-1, q) - \pi(m,l,k+1,k,j,j,q)) -\mu(l)\mu(m) \Big\}\\
&+ \sum_{m=1}^{n-k} \sum_{l=j+1}^{n-k} \sum_{j=k'+1}^k \sum_{q=1}^{\min(k+1, l-m)} |\Gamma^-_{k+1}(m,l,j,q)| \Big\{ \prod_{i=k'+1}^k p_i^{2\binom{k+1}{i+1} - \binom{j}{i+1}} (\pi(m,l,k+1, k+1, j,j, q) \\
&+ \pi(m,l,k, k, j,j, q) - \pi(m,l,k, k+1, j,j, q) - \pi(m,l,k+1,k,j,j,q)) -\mu(l)\mu(m) \Big\} = V_1 + V_2.
\end{align*}

Similarly, we calculate the remaining term in the expansion of the variance \eqref{varexapansion}.
We notice that if $l=m$, then $q=0$ and we have $\Gamma_{k+1}(l,l,j,0) = \Gamma^+_{k+1}(l,l,j,0)$. Hence,  $|\Gamma_{k+1}(l,l,j,0)| = \binom{n-l}{2k+1-j} \binom{2k+1-j}{k} \binom{k}{j-1}$, and
\begin{align*}
&\sum_{\substack{s \neq t \in \binom{[n]}{k+1} \\ \min(s) = \min(t)}} \Cov(Z_sY_s, Z_tY_t) = \sum_{l=1}^{n-k} \sum_{j=k'+1}^k \sum_{(s,t) \in \Gamma_{k+1}(l,l,j,0)} \Cov(Z_sY_s, Z_tY_t) \\
= &\sum_{l=1}^{n-k} \sum_{j=k'+1}^k  \binom{n-l}{2k+1-j} \binom{2k+1-j}{k} \binom{k}{j-1} \Big\{ \prod_{i=k'+1}^k p_i^{2\binom{k+1}{i+1} - \binom{j}{i+1}}  \\
& \left[ \pi(l,l,k+1,k+1,j,0,0) + \pi(l,l,k,k,j-1,0,0) - 2\pi(l,l,k+1,k,j-1,0,0) \right] - \mu(l)^2\Big\}= V_3.
\end{align*}
\end{proof}

\begin{proof}[Proof for Lemma \ref{lemma:lower_bound_var_crit_lexi}]
Fix $k'+ 1 \leq k \leq n-1$ and $p \in (0,1)$, and consider the variance. From Lemma \ref{lemma:var_crit_lexi} we have $\Var\{ T_{k+1} \} = V_1 + V_2 + V_3 + V_4$. 
First we lower bound $V_1$ and $V_2$ by just the negative part of the sum:
\begin{align*}
V_1 \geq  &-2\rho^+(k+1)^2\sum_{l < m}^{n-k} \sum_{j=k'+1}^k \sum_{q=1}^{\min(k+1, m-l)} \binom{n-m}{2k+1-j-q} \binom{2k+1-j-q}{k} \binom{k}{j-1} \binom{m-l+1}{q-1} \\
&\Big\{ \rho^+(j)^{-1} \tau(l,m,q,k,j-1) \left(1-\rho(k+1)-\rho(k) + \rho(k+1)\rho(k)\rho(j-1)^{-1}\right)^{l-1} \\
&+ \rho^+(j)^{-1} \tau(l,m,q,k+1,j)  \left(1-\rho(k+1)-\rho(k) + \rho(k+1)\rho(k)\rho(j)^{-1}\right)^{l-1} \\
&+ (1-\rho(k+1))^{l+m-2} + (1-\rho(k))^{l+m-2} \Big\};\\
V_2 \geq  &-2\rho^+(k+1)^2 \sum_{l < m}^{n-k} \sum_{j=k'+1}^k \sum_{q=1}^{\min(k+1, m-l)} \binom{n-m}{2k+1-j-q} \binom{2k+1-j-q}{k} \binom{k}{j} \binom{m-l+1}{q-1} \\
&\Big\{ \rho^+(j)^{-1} \tau(l,m,q,k,j) \left(1-\rho(k+1)-\rho(k) + \rho(k+1)\rho(k)\rho(j)^{-1}\right)^{l-1}\\
&+ \rho^+(j)^{-1} \tau(l,m,q,k+1,j) \left(1-\rho(k+1)-\rho(k) + \rho(k+1)\rho(k)\rho(j)^{-1}\right)^{l-1}\\
&+ (1-\rho(k+1))^{l+m-2} + (1-\rho(k))^{l+m-2}\Big\}.
\end{align*}

Now using that $\binom{k}{j} + \binom{k}{j-1} = \binom{k+1}{j}$ and $\rho(a) \leq \rho(b)$ for any two integers $b < a$, it is easy to see that $V_1 + V_2 \geq -8\rho^+(k+1)^2R_1 - 8\rho^+(k+1)^2R_2$, where
\begin{align*}
    &R_1 \coloneqq \sum_{l < m}^{n-k} \sum_{j=k'+1}^k \sum_{q=1}^{\min(k+1, m-l)} \binom{n-m}{2k+1-j-q} \binom{2k+1-j-q}{k} \binom{k+1}{j} \binom{m-l+1}{q-1} (1-\rho(k+1))^{l+m-2}; \\
    &R_2 \coloneqq \sum_{l < m}^{n-k} \sum_{j=k'+1}^k \sum_{q=1}^{\min(k+1, m-l)} \binom{n-m}{2k+1-j-q} \binom{2k+1-j-q}{k} \binom{k+1}{j} \binom{m-l+1}{q-1} \\
    &\rho^+(j)^{-1}(1-\rho(k+1))^{m-l}\left(1-\rho(k+1)-\rho(k) + \rho(k+1)\rho(k)\rho(j)^{-1}\right)^{l-1}.
\end{align*}

For $V_3$ we lower bound by terms with $m = k'+1$ and the negative parts of the other terms. Note that $\rho^+(k'+1) = 1 = \rho(k')$ and $\rho(k'+1) < 1$. Now we can proceed:
\begin{align*}
V_3 &\geq \rho^+(k+1)^2\sum_{l=1}^{n-k} \binom{n-l}{2k - k'} \binom{2k - k'}{k} \binom{k}{k'}  \left\{ (1-2\rho(k+1) + \rho(k+1)^2\rho(k'+1)^{-1})^{l-1} -  (1-\rho(k+1))^{2l-2}\right\}\\
&-2\rho^+(k+1)^2\sum_{l=1}^{n-k} \sum_{j=k'+2}^k  \binom{n-l}{2k+1-j} \binom{2k+1-j}{k} \binom{k}{j-1} \\
&\Big\{ \rho^+(j)^{-1}\left(1 -\rho(k) -\rho(k+1) + \rho(k)\rho(k+1)\rho(j-1)^{-1}\right)^{l-1} +  (1-\rho(k+1))^{2l-2}\Big\}\\
&=\rho^+(k+1)^2R_4 - 2\rho^+(k+1)^2R_3;
\end{align*}
here we call the first sum of the lower bound $R_4$ and the second sum $R_3$. For $V_4$ we use the trivial lower bound $V_4 \geq 0$. Hence, we have:
\[\Var(T_{k+1}) \geq \rho^+(k+1)^2(R_4 - 8R_1 - 8R_2 - 2R_3).\]

Let us now upper bound $R_1$:
\begin{align*}
R_1 &\leq \sum_{l < m}^{n-k} \sum_{j=k'+1}^k \sum_{q=1}^{\min(k+1, m-l)} \frac{(n-m)^{2k+1-j-q}}{(2k+1-j-q)!} \frac{(2k+1-j-q)^k}{k!} \frac{(k+1)^j}{j!} \frac{(m-l+1)^{q-1}}{(q-1)!} (1-\rho(k+1))^{l+m-2}\\
&\leq \sum_{l < m}^{n-k} \sum_{q=1}^{\min(k+1, m-l)}(k-k') \frac{n^{2k-k'-q}}{1} \frac{(2k-k'-1)^k}{k!} \frac{(k+1)^k}{(k'+1)!} \frac{n^{q-1}}{1} (1-\rho(k+1))^{l+m-2}\\
&\leq n^{2k-k'-1}(k-k') \frac{(2k-k'-1)^k}{k!} \frac{(k+1)^{k+1}}{(k'+1)!} \sum_{l < m}^{\infty} (1-\rho(k+1))^{l+m-2}\\
&= n^{2k-k'-1}(k-k') \frac{(2k-k'-1)^k}{k!} \frac{(k+1)^{k+1}}{(k'+1)!} \frac{1 - \rho(k+1)}{(2-\rho(k+1))\rho(k+1)^2}.
\end{align*}

Noting that $(1-\rho(k+1))^{m-l}\left(1-\rho(k+1)-\rho(k) + \rho(k+1)\rho(k)\rho(j)^{-1}\right)^{l-1} \leq (1-\rho(k+1))^{m-1}$ (since $\rho(k+1)\rho(j)^{-1} \leq 1$ for all $j \in [k]$), we can bound $R_2$ in an identical way:
\begin{align*}
    R_2 &\leq n^{2k-k'-1}(k-k') \frac{(2k-k'-1)^k}{k!} \frac{(k+1)^{k+1}}{(k'+1)!} \rho^+(k)^{-1} \sum_{l < m}^\infty (1-\rho(k+1))^{m-1} \\
    &= n^{2k-k'-1}(k-k') \frac{(2k-k'-1)^k}{k!} \frac{(k+1)^{k+1}}{(k'+1)!} \frac{1-\rho(k+1)}{\rho^+(k)\rho(k+1)^2}.
\end{align*}

Noting that $\left(1 -\rho(k) -\rho(k+1) + \rho(k)\rho(k+1)\rho(j-1)^{-1}\right)^{l-1}  \leq (1-\rho(k+1))^{l-1}$ and $(1-\rho(k+1))^{2l-2} \leq (1-\rho(k+1))^{l-1}$ we proceed to bound $R_3$:
\begin{align*}
    R_3 &\leq \sum_{l=1}^{n-k} \sum_{j=k'+2}^k  \binom{n-l}{2k+1-j} \binom{2k+1-j}{k} \binom{k}{j-1} (\rho^+(j)^{-1} + 1) (1-\rho(k+1))^{l-1}\\
    &\leq \sum_{l=1}^{n-k} \sum_{j=k'+2}^k  \frac{n^{2k+1-j}}{(2k+1-j)!} \frac{(2k+1-j)^k}{k!} \frac{k^{j-1}}{(j-1)!} (\rho^+(k)^{-1} + 1) (1-\rho(k+1))^{l-1}\\
    &\leq  (k-k'-1)  \frac{n^{2k-k'-1}}{(k+1)!} \frac{(2k-k'-1)^k}{k!} \frac{k^{k}}{(k'+1)!} (\rho^+(k)^{-1} + 1) \sum_{l=1}^{\infty} (1-\rho(k+1))^{l-1}\\
    &=n^{2k-k'-1}  \frac{k-k'-1}{(k+1)!} \frac{(2k-k'-1)^k}{k!} \frac{k^{k}}{(k'+1)!} (\rho^+(k)^{-1} + 1) \rho(k+1)^{-1}.
\end{align*}
To lower bound $R_4$ we just take the $l = 2$ term as all summands are non-negative:
\begin{align*}
R_4 &\geq \binom{n-2}{2k - k'} \binom{2k - k'}{k} \binom{k}{k'}  \left\{  \rho(k+1)^2\rho(k'+1)^{-1} -   \rho(k+1)^2\right\} \\
& \geq \frac{(n-2)^{2k-k'}}{(2k-k')^{2k-k'}} \binom{2k-k'}{k} \binom{k}{k'} \rho(k+1)^2(\rho(k'+1)^{-1}-1).
\end{align*}

Since $R_1$, $R_2$, $R_3$ are all at most of the order $n^{2k-k'-1}$ and $R_4$ is at least of the order
$n^{2k-k'}$, we have that for any fixed $k \geq k' + 1$ and a sequence of probabilities $p_i \in (0,1)$
there exists a constant $C > 0$ independent of $n$ and a natural number $N$ such that for any $n \geq N$:
\[\Var(T_{k+1}) \geq \rho^+(k+1)^2(R_4 - 8R_1 - 8R_2 - 2R_3) \geq Cn^{2k-k'}.\]
\end{proof}

\begin{proof}[Proof of Lemma \ref{theorem:crit_lexi_limiting_cov}]
For $i < j$, write $k'+i = k$ and $k'+j = r$. Note that $\sigma_{k}^2$ and $\sigma_{r}^2$ are sums and so is $\Cov(T_{k+1}, T_{r+1})$. Hence, to get the limit of the ratio we can just look at the ratio of the highest order terms from each sum. From the proof of Lemma \ref{lemma:lower_bound_var_crit_lexi} it is clear that the highest order term for $\sigma_k^2$ is the sum of the terms in $V_3$, where the index of the sum $j$ has the value $k'+1$. That is,
\begin{align*}
    \sigma_{k}^2 &= \sum_{l=1}^{n-k} \binom{n-l}{2k-k'}\binom{2k-k'}{k}\binom{k}{k'}\rho^+(k+1)^2 \\
    &\left[ (1-2\rho(k+1) + \rho(k+1)^2\rho(k'+1)^{-1})^{l-1} - (1-2\rho(k+1) + \rho(k+1)^2)^{l-1} \right] + o(n^{2k-k'}).
\end{align*}

Analogously to the calculations in Lemmas \ref{lemma:var_crit_lexi} and \ref{lemma:lower_bound_var_crit_lexi}, one can perform calculations to find the highest order term of $\Cov(T_{k+1}, T_{r+1})$. Similarly, like in the variance case, it will correspond to subsets $s \in \binom{[n]}{k+1}, t \in \binom{[n]}{r+1}$ such that $\min(s) = \min(t)$ and $|s \cap t| = k'+1$. That is, borrowing the notation from the proof of Lemma \ref{lemma:var_crit_lexi}, we have
\begin{align*}
    &\Cov(T_{k+1}, T_{r+1}) = \sum_{\substack{s \in \binom{[n]}{k+1}, t \in \binom{[n]}{r+1} \\ \min(s) = \min(t) \\ |s \cap t| = k'+1}} \Cov(Z_sY_s, Z_tY_t) + o(n^{k+r-k'}) \\
    &= o(n^{k+r-k'}) + \sum_{l=1}^{n-r} \binom{n-l}{k+r-k'}\binom{k+r-k'}{k}\binom{k}{k'}\rho^+(k+1)\rho^+(r+1) \\
    &\big[ (1-\rho(k+1) -\rho(r+1) + \rho(k+1)\rho(r+1)\rho(k'+1)^{-1})^{l-1}- (1-\rho(k+1) -\rho(r+1) + \rho(k+1)\rho(r+1))^{l-1} \big].
\end{align*}

To deal with those sums, we will use this auxiliary claim:
\begin{claim}\label{claim:aux_crit_lexi_limiting_cov}
For any fixed positive integer $k$ and $x \in (0,1)$ we have:
\[\sum_{l=1}^n \binom{n-l}{k}x^{l-1} = \frac{n^k}{k!}(1-x)^{-1} + o(n^k).\]
\end{claim}

Now we can proceed with the limit, cancelling the binomial coefficients in a similar way to Lemma \ref{lemma:binomial_limit} and then using Claim \ref{claim:aux_crit_lexi_limiting_cov}.

\begin{align*}
    &\lim_{n \to \infty} \Sigma_{i,j} = \\
    &\lim_{n \to \infty} \sum_{l=1}^{n-r} \binom{n-l}{k+r-k'}(k+r-k')!\big[ (1-\rho(k+1) -\rho(r+1) + \rho(k+1)\rho(r+1)\rho(k'+1)^{-1})^{l-1} \\
    &- (1-\rho(k+1) -\rho(r+1) + \rho(k+1)\rho(r+1))^{l-1} \big]\\
    &\left\{\sum_{l=1}^{n-k} \binom{n-l}{2k-k'}(2k-k')! \left[ (1-2\rho(k+1) + \rho(k+1)^2\rho(k'+1)^{-1})^{l-1} - (1-2\rho(k+1) + \rho(k+1)^2)^{l-1} \right]\right\}^{-\frac{1}{2}} \\
    &\left\{\sum_{l=1}^{n-r} \binom{n-l}{2r-k'}(2r-k')! \left[ (1-2\rho(r+1) + \rho(r+1)^2\rho(k'+1)^{-1})^{l-1} - (1-2\rho(r+1) + \rho(r+1)^2)^{l-1} \right]\right\}^{-\frac{1}{2}}= \\
    &\lim_{n \to \infty} n^{k+r-k'}\big[ (\rho(k+1) +\rho(r+1) - \rho(k+1)\rho(r+1)\rho(k'+1)^{-1})^{-1} - (\rho(k+1) +\rho(r+1) - \rho(k+1)\rho(r+1))^{-1} \big]\\
    &\left\{n^{2k-k'} \left[ (2\rho(k+1) - \rho(k+1)^2\rho(k'+1)^{-1})^{-1} - (2\rho(k+1) - \rho(k+1)^2)^{-1} \right]\right\}^{-\frac{1}{2}} \\
    &\left\{n^{2r-k'} \left[ (2\rho(r+1) - \rho(r+1)^2\rho(k'+1)^{-1})^{-1} - (2\rho(r+1) - \rho(r+1)^2)^{-1} \right]\right\}^{-\frac{1}{2}} \\
    &= \frac{\rho(k+1)\rho(r+1) \sqrt{(2 - \rho(k+1))(2 - \rho(r+1))(2 - \rho(k+1)\rho(k'+1)^{-1})(2 - \rho(r+1)\rho(k'+1)^{-1})}}{(\rho(k+1) +\rho(r+1) - \rho(k+1)\rho(r+1)\rho(k'+1)^{-1})(\rho(k+1) +\rho(r+1) - \rho(k+1)\rho(r+1))}.
\end{align*}
\end{proof}

\newpage
\section{A multivariate CLT for dissociated sums}\label{section:appendix_clt}

Let $n$ and $d$ be positive integers. For each $i \in [d] =:\{1, 2, \ldots, d\}$, we fix an index set $\I_i \subset [n] \times \set{i}$ and consider the union of disjoint sets $\I =: \bigcup_{i \in [d]} \I_i.$ Let $X_s, s = (k,i) \in \I$ be a collection of real centered random variables which are defined on the same probability space. For each $i \in [d]$ form the sum
\[
W_i \coloneqq \sum_{s \in \I_i} X_s.
\] 
Our interest here is in the resulting random vector $W = (W_1,\ldots,W_d) \in \RR^d$. The following notion is a natural multivariate generalisation of the dissociated sum that has been studied in \cite{temvcinas2021multivariate, barbour1989central, moginley1975dissociated}.

\begin{definition}\label{def:disrand}
We call $W$ a {\bf vector of dissociated sums} if for each $s \in \I$ and $j \in [d]$ there exists a {\bf dependency neighbourhood} $\D_j(s) \subset \I_j$ satisfying three criteria:
\begin{enumerate}
	\item the difference $\left(W_j - \sum_{u \in \D_j(s)} X_u\right)$ is independent of $X_s$; 
	\item for each $t \in \I$, the quantity $\left(W_j - \sum_{u \in \D_j(s)} X_u - \sum_{v \in \D_j(t) \setminus \D_j(s)} X_v\right)$ is independent of the pair $(X_s,X_t)$;
	\item $X_s$ and $X_t$ are independent if $t \not\in \bigcup_j \D_j(s)$. 
\end{enumerate}
\end{definition}

Let $W$ be a vector of dissociated sums as defined above.  For each $s\in \I$, as $\mathbb{I}_j \cap \mathbb{I}_\ell = \emptyset$ for $j \ne \ell$ by construction, the sets  $\D_j(s), j \in [d]$ are disjoint (although for $s \ne t$, the sets $\D_j(s)$ and $\D_j(t)$ may not be disjoint). {$\D_j(s)$ is the dependency neighbourhood of $X_s$ in the $j$-th component of the vector $W$.} We write $\D(s) = \bigcup_{j \in [d]}\D_j(s)$ for the disjoint union of these of dependency neighbourhoods.

\begin{theorem}[Theorems 2 and 9 in \cite{temvcinas2021multivariate}]\label{theorem:mvn_dissociated_decomp_approx}
Let $h: \RR^d \to \RR$ be any three times continuously differentiable function whose third partial derivatives are Lipschitz continuous and bounded. Consider a standard $d$-dimensional Gaussian vector $Z \sim \mvn{0}{\text{\rm Id}_{d \times d}}$. Assume that for all $s \in \I$, we have $\EEbracket{X_s} = 0$ and $\EE \abs {X_s^3} < \infty.$ Then, for any vector of dissociated sums $W \in \RR^d$ with a positive semi-definite covariance matrix $\Sigma$:
\begin{enumerate}
    \item \[
\abs{\EE h(W) - \EE h(\Sigma^{\frac{1}{2}}Z)} \leq B_{\ref{theorem:mvn_dissociated_decomp_approx}} \abs{h}_3,
\]
    \item \[
\sup _{A \in \mathcal{K}}|\PP(W \in A)-\PP(\Sigma^{\frac{1}{2}}Z \in A)| \leq 2^{\frac{7}{2}} 3^{-\frac{3}{4}}d^{\frac{3}{16}}B_{\ref{theorem:mvn_dissociated_decomp_approx}}^{\frac{1}{4}},
\]
\end{enumerate}

where $B_{\ref{theorem:mvn_dissociated_decomp_approx}} = B_{\ref{theorem:mvn_dissociated_decomp_approx}.1} + B_{\ref{theorem:mvn_dissociated_decomp_approx}.2}$ is the sum given by
\begin{align*}
B_{\ref{theorem:mvn_dissociated_decomp_approx}.1} &\coloneqq 
\frac{1}{3}\sum_{s \in \I} \sum_{t, u \in \D(s)} \hspace{-.4em} \left(\frac{1}{2}\EE\abs{X_sX_tX_u}+\EE\abs{X_sX_t} \EE\abs{X_u}\right) \\
B_{\ref{theorem:mvn_dissociated_decomp_approx}.2} &\coloneqq  \frac{1}{3} \sum_{s \in \I} \sum_{t \in \D(s)} \sum_{v \in \D(t) \setminus \D(s)} \hspace{-.4em} \left(\EE\abs{X_sX_tX_v} + \EE\abs{X_sX_t} \EE\abs{X_v}\right). 
\end{align*}
\end{theorem}

\begin{corollary}[Corollaries 8 and 11 in \cite{temvcinas2021multivariate}]\label{corollary:mvn_dissociated_decomp}
Let the assumptions of Theorem \ref{theorem:mvn_dissociated_decomp_approx} hold. If we further require that the random variables $X_{s}$ and $X_{t}$ are independent whenever $t$ lies in $\I \setminus \D(s)$, then for any vector of dissociated sums $W \in \RR^d$ with a positive semi-definite covariance matrix $\Sigma$,
\begin{enumerate}
    \item \[ \abs{\EE h(W) - \EE h(\Sigma^{\frac{1}{2}}Z)} \leq B_{\ref{corollary:mvn_dissociated_decomp}} \abs{h}_3; \]
    \item  \[\sup _{A \in \mathcal{K}}|\PP(W \in A)-\PP(\Sigma^{\frac{1}{2}}Z \in A)| \leq 2^{\frac{7}{2}} 3^{-\frac{3}{4}}d^{\frac{3}{16}}B_{\ref{corollary:mvn_dissociated_decomp}}^{\frac{1}{4}},
    \]
\end{enumerate}

where $B_{\ref{corollary:mvn_dissociated_decomp}}$ is a sum over $(i,j,k) \in [d]^3$ of the form:
\[
B_{\ref{corollary:mvn_dissociated_decomp}} := \frac{1}{3} \sum_{(i,j,k)} \abs{\I_i}  \alpha_{ij}  \left(\frac{3\alpha_{ik}}{2} + 2\alpha_{jk}\right)  \beta_{ijk}. 
\]
Here $\alpha_{ij}$ is the largest  value attained by $\abs{\D_j(s)}$ over $s \in \I_i$, and 
\[
\beta_{ijk} = \max_{s,t,u}\Big(\EE\abs{X_sX_tX_u},\EE\abs{X_sX_t}\EE\abs{X_u}\Big)
\] as $(s,t,u)$ range over $\I_i \times \I_j \times \I_k$.
\end{corollary}

The following lemma is useful when we deal with Bernoulli variables.
\begin{lemma}[Lemma 9 in \cite{temvcinas2021multivariate}]\label{lemma:unified_moments}
Let $\xi_1, \xi_2, \xi_3$ be Bernoulli random variables with expected values $\mu_1, \mu_2, \mu_3$ respectively. Let $c_1, c_2, c_3 > 0$ be any deterministic quantities. Consider variables $X_i \coloneqq c_i (\xi_i - \mu_i)$ for $i = 1,2,3$. Then we have
\begin{align*}
    &\EE \abs{X_1 X_2 X_3} \leq c_1 c_2 c_3 \left\{ \mu_1\mu_2(1-\mu_1)(1-\mu_2) \right\}^{\frac{1}{2}}; \\
    &\EE \abs{X_{1} X_{2}} \EE \abs{X_{3}} \leq c_1 c_2 c_3 \left\{ \mu_1\mu_2(1-\mu_1)(1-\mu_2) \right\}^{\frac{1}{2}}.
\end{align*}
\end{lemma}

\newpage
\section{Simulation Study Results}\label{section:appendix_simulation}

\begin{table}[!h]
\begin{center}
\begin{tabular}{|c|c|c|c|c|}
\hline $\epsilon_1$ & $\epsilon_2$ & centered triangles & triangles  & critical \\
\hline 0.00818181818181818 & 0.462727272727272 & 96 & 100 & 99 \\
\hline 0.0163636363636363 & 0.425454545454545 & 94 & 100 & 100 \\
\hline 0.0245454545454545 & 0.388181818181818 & 88 & 100 & 100 \\
\hline 0.0327272727272727 & 0.35090909090909 & 93 & 100 & 99 \\
\hline 0.0409090909090909 & 0.313636363636363 & 94 & 100 & 100 \\
\hline 0.049090909090909 & 0.276363636363636 & 96 & 100 & 99 \\
\hline 0.0572727272727272 & 0.239090909090909 & 93 & 100 & 99 \\
\hline 0.0654545454545454 & 0.201818181818181 & 95 & 100 & 100 \\
\hline 0.0736363636363636 & 0.164545454545454 & 91 & 100 & 100 \\
\hline 0.0818181818181818 & 0.127272727272727 & 95 & 100 & 81 \\
\hline 0.0825619834710743 & 0.123884297520661 & 93 & 100 & 69 \\
\hline 0.0833057851239669 & 0.120495867768595 & 91 & 100 & 67 \\
\hline 0.0840495867768594 & 0.117107438016528 & 93 & 100 & 50 \\
\hline 0.084793388429752 & 0.113719008264462 & 94 & 100 & 40 \\
\hline 0.0855371900826446 & 0.110330578512396 & 93 & 100 & 29 \\
\hline 0.0862809917355371 & 0.10694214876033 & 96 & 100 & 23 \\
\hline 0.0870247933884297 & 0.103553719008264 & 92 & 100 & 18 \\
\hline 0.0877685950413223 & 0.100165289256198 & 97 & 100 & 12 \\
\hline 0.0885123966942148 & 0.0967768595041322 & 97 & 100 & 14 \\
\hline 0.0892561983471074 & 0.0933884297520661 & 89 & 100 & 5 \\
\hline 0.09 & 0.09 & 98 & 100 & 2 \\
\hline
\end{tabular}
\caption{\label{table:dim3} Goodness of fit testing results of the tetrahedron model.}
\end{center}
\end{table}

\begin{table}[!h]
\begin{center}
\begin{tabular}{|c|c|c|c|c|}
\hline $\epsilon_1$ & $\epsilon_2$ & centered triangles & triangles  & critical \\
\hline 0.00818181818181818 & 0.608181818181818 & 96 & 100 & 99 \\
\hline 0.0163636363636363 & 0.556363636363636 & 95 & 100 & 98 \\
\hline 0.0245454545454545 & 0.504545454545454 & 95 & 100 & 95 \\
\hline 0.0327272727272727 & 0.452727272727272 & 92 & 100 & 100 \\
\hline 0.0409090909090909 & 0.40090909090909 & 90 & 100 & 98 \\
\hline 0.049090909090909 & 0.349090909090909 & 96 & 100 & 100 \\
\hline 0.0572727272727272 & 0.297272727272727 & 97 & 100 & 98 \\
\hline 0.0654545454545454 & 0.245454545454545 & 93 & 100 & 93 \\
\hline 0.06681818182 & 0.2368181819 & 91 & 100 & 86 \\
\hline 0.06818181818 & 0.2281818182 & 94 & 100 & 89 \\
\hline 0.069545454545 & 0.21954545455 & 95 & 100 & 81 \\
\hline 0.07090909091 & 0.2109090909 & 96 & 100 & 72 \\
\hline 0.072272727275 & 0.20227272725 & 92 & 100 & 65 \\
\hline 0.0736363636363636 & 0.193636363636363 & 98 & 100 & 67 \\
\hline 0.0751239669421487 & 0.184214876033057 & 96 & 100 & 63 \\
\hline 0.0766115702479338 & 0.174793388429752 & 95 & 100 & 55 \\
\hline 0.078099173553719 & 0.165371900826446 & 98 & 100 & 53 \\
\hline 0.0795867768595041 & 0.15595041322314 & 97 & 100 & 45 \\
\hline 0.0810743801652892 & 0.146528925619834 & 94 & 100 & 38 \\
\hline 0.0825619834710743 & 0.137107438016528 & 96 & 100 & 28 \\
\hline 0.0840495867768594 & 0.127685950413223 & 95 & 100 & 15 \\
\hline 0.0855371900826446 & 0.118264462809917 & 96 & 100 & 12  \\
\hline 0.0870247933884297 & 0.108842975206611 & 87 & 100 & 11 \\
\hline 0.09 & 0.0994214876033057 & 97 & 100 & 9 \\
\hline
\end{tabular}
\caption{\label{table:dim2} Goodness of fit testing results of the triangle model.}
\end{center}
\end{table}

\begin{table}[!h]
\begin{center}
\begin{tabular}{|c|c|c|c|c|}
\hline $\epsilon_1$ & $\epsilon_2$ & centered triangles & triangles & critical \\
\hline 0 & 1.732050808 & 100 & 100 & 100 \\
\hline 0.02344631004 & 1.673018508 & 99 & 100 & 99 \\
\hline 0.04689262009 & 1.613986208 & 90 & 100 & 100 \\
\hline 0.07033893013 & 1.554953908 & 97 & 100 & 99 \\
\hline 0.09378524018 & 1.495921608 & 96 & 100 & 99 \\
\hline 0.1172315502 & 1.436889308 & 97 & 100 & 99 \\
\hline 0.1406778603 & 1.377857009 & 94 & 100 & 100 \\
\hline 0.1641241703 & 1.318824709 & 96 & 100 & 100 \\
\hline 0.1875704804 & 1.259792409 & 99 & 100 & 100 \\
\hline 0.2110167904 & 1.200760109 & 91 & 100 & 100 \\
\hline 0.2344631004 & 1.141727809 & 95 & 100 & 99 \\
\hline 0.2579094105 & 1.082695509 & 77 & 100 & 97 \\
\hline 0.2813557205 & 1.023663209 & 43 & 100 & 98 \\
\hline 0.2857519 & 1.01259465 & 39 & 100 & 99 \\
\hline 0.29014809 & 1.0015261 & 26 & 100 & 97 \\
\hline 0.29454427 & 0.99045754 & 17 & 100 & 93 \\
\hline 0.29894045 & 0.97938898 & 12 & 100 & 88 \\
\hline 0.30333664 & 0.96832043 & 5 & 100 & 75 \\
\hline 0.3048020306 & 0.9646309096 & 0 & 100 & 73 \\
\hline 0.30773282 & 0.95725187 & 1 & 100 & 81 \\
\hline 0.312129 & 0.94618332 & 0 & 98 & 62 \\
\hline 0.31652519 & 0.93511476 & 0 & 96 & 51 \\
\hline 0.32092137 & 0.9240462 & 0 & 90 & 33 \\
\hline 0.32531755 & 0.91297765 & 0 & 75 & 29 \\
\hline 0.3282483406 & 0.9055986098 & 3 & 74 & 12 \\
\hline 0.32971373 & 0.90190909 & 0 & 53 & 9 \\
\hline 0.33410992 & 0.89084053 & 0 & 33 & 4 \\
\hline 0.3385061 & 0.87977198 & 0 & 16 & 3 \\
\hline 0.34290228 & 0.86870342 & 0 & 4 & 1 \\
\hline 0.34729847 & 0.85763487 & 0 & 3 & 1 \\
\hline 0.3516946507 & 0.84656631 & 0 & 3 & 0 \\
\hline 0.3751409607 & 0.7875340101 & 0 & 0 & 0 \\
\hline 0.3985872707 & 0.7285017103 & 0 & 0 & 0 \\
\hline 0.4220335808 & 0.6694694104 & 0 & 0 & 0 \\
\hline 0.4454798908	& 0.6104371106 & 0 & 0 & 0 \\
\hline 0.4689262009 &	0.5514048108 & 0 & 0 & 0 \\
\hline 0.4923725109 &	0.4923725109 & 0 & 0 & 0 \\
\hline
\end{tabular}
\caption{\label{table:dim1} Goodness of fit testing results of the edge model.}
\end{center}
\end{table}

\end{appendix}
\end{document}